\documentclass{imanum}
\usepackage{graphicx}
\usepackage{multirow,bigstrut}
\usepackage{amsmath,amsfonts,amssymb}
\usepackage{mathrsfs}
\usepackage{color}
\usepackage{upgreek}
\usepackage{bm}
\usepackage{verbatim}

\jno{drnxxx}
\begin{document}

\title{Error estimates of the Crank-Nicolson Galerkin method for
the time-dependent Maxwell-Schr\"{o}dinger equations under the Lorentz gauge}
% Short title for running heads:
\shorttitle{Crank-Nicolson Galerkin method for Maxwell-Schr\"{o}dinger Equations}

\author{%
{\sc
Chupeng Ma\thanks{Email: machupeng@lsec.cc.ac.cn},
Liqun Cao\thanks{Corresponding author. Email: clq@lsec.cc.ac.cn}
} \\[2pt]
Institute of Computational Mathematics and Scientific/Engineering Computing,\\
      Academy of Mathematics and Systems Science, \\
      Chinese Academy of  Sciences, Beijing 100190, China \\[6pt]
{\sc and}\\[6pt]
{\sc Yanping Lin}\thanks{Email: malin@polyu.edu.hk}\\[2pt]
Department of Applied Mathematics,\\
The Hong Kong Polytechnic University, Hung Hom, Kowloon, Hong Kong, China
}
% Short list of authors for running heads:
\shortauthorlist{C. P. Ma \emph{et al.}}

\maketitle

\begin{abstract}
% Body of abstract:
{In this paper we study the numerical method and the convergence for solving the time-dependent Maxwell-Schr\"{o}dinger equations under the Lorentz gauge. An alternating Crank-Nicolson finite element method for solving the problem is presented and
the optimal error estimate for the numerical algorithm is obtained by a mathematical inductive method.
Numerical examples are then carried out to confirm the theoretical results.}
% Keywords:
{error estimates; Crank-Nicolson; Galerkin method;
Maxwell-Schr\"{o}dinger.}
\end{abstract}

\section{Introduction}\label{sec-1}
Light-matter interaction at nanoscale is a central topic in the study of optical properties of nanophotonic systems, for example, metallic nanostructures and quantum dots.
In view of practical numerical simulation,  a semiclassic model is often used for modelling light-matter interaction. The basic idea is to use the classical Maxwell's equations for the electromagnetic field and the Schr\"{o}dinger equation for the matter. In this paper, we study the following Maxwell-Schr\"{o}dinger coupled system, which describes the interaction between an
electron and  its  self-consistent generated and external electromagnetic fields.
\begin{equation}\label{eq:1-1}
\left\{
\begin{array}{@{}l@{}}
{\displaystyle  \mathrm{i}\hbar\frac{\partial \psi(\mathbf{x},t)}{\partial t}=
\left\lbrace\frac{1}{2m}\left[\mathrm{i}\hbar\nabla +q\mathbf{A}(\mathbf{x},t)\right]^{2}
 + q \phi(\mathbf{x},t)+V_{0} \right\rbrace\psi(\mathbf{x},t),}\\[2mm]
{\displaystyle  \qquad \qquad \qquad \qquad \qquad \qquad \qquad \qquad \qquad \qquad(\mathbf{x},t)\in\Omega\times(0,T),}\\[2mm]
 {\displaystyle -\frac{\partial}{\partial t}\nabla\cdot\big(\epsilon\mathbf{A}(\mathbf{x},t)\big)
 -\nabla\cdot\big(\epsilon\nabla\phi(\mathbf{x},t)\big) =q |\psi(\mathbf{x},t)|^{2}, \,\,
 (\mathbf{x},t)\in\Omega\times(0,T),} \\[2mm]
{\displaystyle \epsilon\frac{\partial ^{2}\mathbf{A}(\mathbf{x},t)}{\partial t^{2}}+\nabla\times
\big({\mu}^{-1}\nabla\times \mathbf{A}(\mathbf{x},t)\big)
+\epsilon \frac{\partial (\nabla \phi(\mathbf{x},t))}{\partial t} =\mathbf{J}_{q}(\mathbf{x},t),}\\[2mm]
 {\displaystyle  \qquad \qquad \qquad \qquad \qquad \qquad \qquad \qquad \qquad \qquad
(\mathbf{x},t)\in \Omega\times(0,T),}\\[2mm]
{\displaystyle \mathbf{J}_q=-\frac{\mathrm{i}q\hbar}{2m}\big(\psi^{\ast}\nabla{\psi}-\psi\nabla{\psi}^{\ast}\big)-\frac{\vert q\vert^{2}}{m}\vert\psi\vert^{2}\mathbf{A},} \\[2mm]
{\displaystyle \psi, \phi, \mathbf{A} \,\, \mathrm{subject \ to \ the \ appropriate \ initial \ and\ boundary \ conditions}, }
\end{array}
\right.
\end{equation}
where $ \Omega\subset \mathbb{R}^d  $, $ d\geq 2 $ is a bounded Lipschitz polygonal convex domain in $\mathbb{R}^{2}$
(or a bounded Lipschitz polyhedron convex domain
in $\mathbb{R}^{3}$),
$ \psi^{\ast} $ denotes the complex conjugate of $ \psi $,
$\epsilon$ and $\mu$ respectively denote the electric permittivity and the magnetic permeability of the material and $V_0$ is the constant potential energy.

It is well-known that the solutions of the Maxwell-Schr\"{o}dinger equations (\ref{eq:1-1}) lack uniqueness. In fact,
for any function $\chi : \Omega\times (0,T)\rightarrow \mathbb{R}$, if $(\psi, \phi, \mathbf{A}) $ is a solution of (\ref{eq:1-1}), then $(\exp(i\chi)\psi, \phi-\partial_{t}\chi, \mathbf{A}+\nabla \chi) $ is also a solution of (\ref{eq:1-1}). To obtain mathematically well-posed equations, some extra constraint, commonly known as gauge choice, is often enforced on the solutions of the Maxwell-Schr\"{o}dinger equations. The most common gauges are listed below.

(i) The Lorentz gauge
\begin{equation*}
{\displaystyle\nabla\cdot\mathbf{A}+\frac{\partial \phi}{\partial t} =0}.
\end{equation*}

(ii) The Coulomb gauge
 \begin{equation*}
{\displaystyle \nabla\cdot\mathbf{A}=0}.
 \end{equation*}

(iii) The temporal gauge
 \begin{equation*}
{\displaystyle\phi=0}.
\end{equation*}

For simplicity, we employ the atomic units,  i.e. $\hbar=m=q=1 $,  and assume that $\epsilon=\mu=1$ without loss of generality.

In this paper, we consider the time-dependent Maxwell-Schr\"{o}dinger equations under the Lorentz gauge as follows:
\begin{equation}\label{eq:1-2}
\left\{
\begin{array}{@{}l@{}}
{\displaystyle  -\mathrm{i}\frac{\partial \psi}{\partial t}+
\frac{1}{2}\left(\mathrm{i}\nabla +\mathbf{A}\right)^{2}\psi
 + V_{0}\psi + \phi\psi= 0 ,\,\, (\mathbf{x},t)\in
\Omega\times(0,T),}\\[2mm]
{\displaystyle \frac{\partial ^{2}\mathbf{A}}{\partial t^{2}}+\nabla\times
(\nabla\times \mathbf{A}) - \nabla(\nabla \cdot \mathbf{A})+\frac{\mathrm{i}}{2}\big(\psi^{*}\nabla{\psi}-\psi\nabla{\psi}^{*}\big) }\\[2mm]
{\displaystyle \qquad\qquad+\vert\psi\vert^{2}\mathbf{A}=0,
\,\,\quad (\mathbf{x},t)\in \Omega\times(0,T),}\\[2mm]
{\displaystyle  \frac{\partial ^{2}\phi}{\partial t^{2}}-\Delta \phi = \vert\psi\vert^{2},\,\, (\mathbf{x},t)\in \Omega\times(0,T) }.
\end{array}
\right.
\end{equation}

The boundary conditions are
\begin{equation}\label{eq:1-3}
{\displaystyle \psi(\mathbf{x},t)=0,\,\,\phi(\mathbf{x},t)=0,\,\, \mathbf{A}(\mathbf{x},t)\times\mathbf{n}=0, \,\, \nabla \cdot \mathbf{A}(\mathbf{x},t) = 0, \quad (\mathbf{x},t)\in
\partial \Omega\times(0,T),}
\end{equation}
and the initial conditions are
\begin{equation}\label{eq:1-4}
\begin{array}{@{}l@{}}
{\displaystyle \psi(\mathbf{x},0) = \psi_0(\mathbf{x}),\quad\phi(\mathbf{x},0) = \phi_0(\mathbf{x}),\quad\phi_{t}(\mathbf{x},0) = \phi_1(\mathbf{x}),}\\[2mm]
{\displaystyle \qquad\qquad \mathbf{A}(\mathbf{x},0)=\mathbf{A}_{0}(\mathbf{x}),\quad\mathbf{A}_{t}(\mathbf{x},0)=\mathbf{A}_{1}(\mathbf{x}),}
\end{array}
\end{equation}
where $ \phi_{t}$ and $ \mathbf{A}_{t} $ denote the derivative of $\phi$ and  $ \mathbf{A} $ with respect to the time $ t $ respectively,
$ \mathbf{n}=(n_1, n_2, n_3) $ is the outward unit normal to the boundary $ \partial \Omega $.

The gauge choice and the equations (\ref{eq:1-2}) impose the following constraints on the initial datas:
\begin{equation}\label{eq:1-5}
{\displaystyle \nabla\cdot\mathbf{A}_{0} + \phi_{1} =0, \quad \nabla\cdot\mathbf{A}_{1} +\Delta\phi_{0} + \vert\psi_{0}\vert^{2} =0.}
\end{equation}

\begin{remark}\label{rem1-1}
The boundary condition $
\psi(x)=0 $ on $\partial \Omega $ implies that the particle is confined in a whole domain $\Omega$. The boundary condition
$  \mathbf{A}(\mathbf{x},t)\times\mathbf{n}=0 $ and $\phi(\mathbf{x},t) = 0$ on $\partial \Omega$  are direct results of  the perfect conductive boundary (PEC).
The boundary condition $\nabla \cdot \mathbf{A}=0$ on $\partial \Omega$ can be deduced from the boundary condition of
$\,\phi$ and the gauge choice. As for the determination of the boundary conditions for the vector potential $\mathbf{A}$ and the scalar potential $\phi$, we refer to \cite{Weng}.
\end{remark}

The local and global well-posedness of solutions on all of $\mathbb{R}^{3}$ for the time-dependent Maxwell-Schr\"{o}dinger equations (\ref{eq:1-1}) have been studied in, for example, \cite{Gin}, \cite{Guo}, \cite{Nak}, \cite{Nak-1}, \cite{Nak-2}, \cite{Shi} and \cite{Wa}. To the best of our knowledge, the existence  and uniqueness of the solution for the Maxwell-Schr\"{o}dinger equations
in a bounded domain seem to be open.

There are a number of results for the numerical methods for the coupled Maxwell-Schr\"{o}dinger equations.
We recall some interesting studies.
\cite{Sui} proposed the finite-difference time-domain (FDTD) method to
solve the Maxwell-Schr\"{o}dinger equations in the simulation of electron tunneling problem.
\cite{Pi} studied  a carbon nanotube between two metallic electrodes by solving the Maxwell-Schr\"{o}dinger equations with the transmission line matrix (TLM)/FDTD hybrid method. \cite{Ah-1} used the FDTD method for the Maxwell-Schr\"{o}dinger system
to simulate plasmonics nanodevices.  However, in the numerical studies listed above,
the EM fields are all described by the Maxwell's equations
involving electric fields $\mathbf{E}$ and magnetic fields $\mathbf{H}$, instead of the $\mathbf{A}$-$\phi$ formulation. Recently, \cite{Ryu} applied the FDTD scheme to discretize the Maxwell-Schr\"{o}dinger equations (\ref{eq:1-1}) under the Lorentz gauge and to simulate the interaction between a single electron in an artificial atom and an incoming electromagnetic field. For more results on this topic, we refer to \cite{Oh}, \cite{Sa}, \cite{Tur} and the references therein.

There are few results on the finite element method (FEM) and its convergence analysis of the Maxwell-Schr\"{o}dinger equations (\ref{eq:1-1}). In this paper we will present an alternating  Crank-Nicolson finite element method for solving the problem
(\ref{eq:1-2})-(\ref{eq:1-4}), i.e. the finite element method in space and the Crank-Nicolson scheme in time.
Then we will derive the optimal error estimates for the proposed method.
Our work is motivated by \cite{Gao-1} in which Gao and his collaborators proposed a linearized Crank-Nicolson Galerkin method for the time-dependent Ginzburg-Landau equations and derived an optimal error estimate via a mathematical inductive method under the assumption that  $h$ and $\Delta t$ are sufficiently small. Here $h$ and $\Delta t$ are the spatial mesh size and the time step, respectively. Compared to the time-dependent Ginzburg-Landau model, the error analysis of numerical schemes for the time-dependent Maxwell-Schr\"{o}dinger system is much more difficult. The main difficulties and tricky parts in this paper are the estimates of the current term $\mathbf{J}_{q}$ and the error analysis for the wave function $\psi$.
In particular we derive the energy-norm error estimates for the Schr\"{o}dinger's equation.

The remainder of this paper is organized as follows.
In section~\ref{sec-2}, we introduce some notation and propose a decoupled alternating Crank-Nicolson scheme with the Galerkin finite element approximation for the Maxwell-Schr\"{o}dinger equations (\ref{eq:1-2})-(\ref{eq:1-4}). The proof of the main theorem (see Theorem~\ref{thm2-1}) in this paper will be given in section ~\ref{sec-3}. Finally, some numerical tests are carried out to validate the theoretical results in this paper.

Throughout the paper the Einstein summation convention on repeated
indices is adopted. By $C$ we shall denote a positive constant
independent of the mesh size $ h $ and the time step $ \Delta t $ without distinction.

\section{An alternating Crank-Nicolson Galerkin finite element scheme}\label{sec-2}

In this section, we present a numerical scheme for the Maxwell-Schr\"{o}dinger equations (\ref{eq:1-2})-(\ref{eq:1-4}) using Galerkin
finite element method in space and a decoupled alternating Crank-Nicolson scheme in time.
To start with, here and afterwards, we assume that $\Omega $ is a bounded Lipschitz polygonal convex domain in $\mathbb{R}^{2}$
(or a bounded Lipschitz polyhedron convex domain in $\mathbb{R}^{3}$).
We introduce the following notation. Let $W^{s,p}(\Omega)$ denote the conventional Sobolev spaces of the real-valued functions.
As usual, $W^{s,2}(\Omega)$ and  $W^{s,2}_{0}(\Omega)$ are denoted by  $H^{s}(\Omega)$ and $H^{s}_{0}(\Omega)$ respectively.
We use $\mathcal{W}^{s,p}(\Omega)=\{u+\mathrm{i}v\,|\, u,v \in W^{s,p}(\Omega)\} $
and $\mathcal{H}^{s}(\Omega)=\{u+\mathrm{i}v\,|\, u,v \in H^{s}(\Omega)\}$ with calligraphic letters for Sobolev spaces of the complex-valued functions, respectively.
Furthermore, let $\mathbf{W}^{s,p}(\Omega) =[W^{s,p}(\Omega)]^{d} $ and $\mathbf{H}^{s}(\Omega)=[H^{s}(\Omega)]^{d}$ with bold faced letters be Sobolev spaces of
the vector-valued functions with $d$ components ($d$=2,\,3). $ L^{2}$ inner-products in $H^{s}(\Omega) $, $\mathcal{H}^{s}(\Omega)$
and $\mathbf{H}^{s}(\Omega)$  are denoted by $(\cdot,\cdot )$ without ambiguity.

In particular, we introduce the following subspace of $\mathbf{H}^{1}(\Omega)$:
\begin{equation*}
\mathbf{H}^{1}_{t}(\Omega)=\{\mathbf{A}\,|\,\mathbf{A}\in\mathbf{H}^{1}(\Omega),\,\,\mathbf{A}\times\mathbf{n}=0 \,\, \,{\rm on}\, \,\,\partial\Omega\}
\end{equation*}
The semi-norm on $\mathbf{H}^{1}_{t}(\Omega)$ is defined by
\begin{equation*}
\Vert \mathbf{u} \Vert_{ \mathbf{H}^{1}_{t}(\Omega)} : = \left[\Vert\nabla \cdot \mathbf{u}\Vert^{2}_{L^{2}(\Omega)} + \Vert\nabla \times \mathbf{u}\Vert^{2}_{\mathbf{L}^{2}(\Omega)
}\right]^{\frac{1}{2}},
\end{equation*}
which is equivalent to the standard $\mathbf{H}^{1}(\Omega)$-norm $\Vert \mathbf{u} \Vert_{ \mathbf{H}^{1}(\Omega)}$, see \cite{Gir}.

The weak formulation of the Maxwell-Schr\"{o}dinger system (\ref{eq:1-2})-(\ref{eq:1-4}) can be specified as follows:
find $(\psi,\mathbf{A},\phi)\in \mathcal{H}_{0}^{1}(\Omega)\times \mathbf{H}_{t}^{1}(\Omega)\times H_{0}^{1}(\Omega)$ such that $\forall t\in(0,T)$,
\begin{equation}\label{eq:2-1}
\left\{
\begin{array}{@{}l@{}}
{\displaystyle
 -\mathrm{i}(\frac{\partial \psi}{\partial t},\varphi)+
\frac{1}{2}(\left(\mathrm{i}\nabla +\mathbf{A}\right)\psi,\left(\mathrm{i}\nabla +\mathbf{A}\right)\varphi)
 + (V_{0}\psi,\varphi)+ (\phi\psi,\varphi) = 0 ,\quad \forall \varphi\in\mathcal{H}_{0}^{1}(\Omega),}\\[2mm]
{\displaystyle ( \frac{\partial ^{2}\mathbf{A}}{\partial t^{2}},\mathbf{v})+(\nabla\times
 \mathbf{A},\nabla\times\mathbf{v}) + (\nabla\cdot
 \mathbf{A},\nabla\cdot\mathbf{v})+(\frac{\mathrm{i}}{2}\big(\psi^{*}\nabla{\psi}-\psi\nabla{\psi}^{*}\big),\mathbf{v}) }\\[2mm]
{\displaystyle\quad \quad  \quad + \,( |\psi|^{2}\mathbf{A},\mathbf{v})=0 ,\qquad \qquad \qquad\qquad\qquad\forall\mathbf{v}\in\mathbf{H}^{1}_{t}(\Omega),}\\[2mm]
{\displaystyle  (\frac{\partial^{2}\phi}{\partial t^{2}},\eta) + (\nabla \phi,\nabla\eta) = (\vert\psi\vert^{2},\eta), \quad\quad \forall \eta\in H_{0}^{1}(\Omega)}
\end{array}
\right.
\end{equation}
with the initial conditions $\psi_{0} \in \mathcal{H}^{1}_{0}(\Omega)$, $\mathbf{A}_{0} \in \mathbf{H}^{1}_{t}(\Omega)$, $\phi_{0}\in H_{0}^{1}(\Omega)$, $\mathbf{A}_{t}(\cdot , 0) \in \mathbf{L}^{2}(\Omega)$ and $\phi_{t}(\cdot,0)\in L^{2}(\Omega)$.

Let M be a positive integer and let $\Delta t =T / M$  be the time step. For any k=1,2,$\cdots, M$, we introduce the following notation:
\begin{equation}\label{eq:2-2}
\begin{array}{@{}l@{}}
{\displaystyle \partial U^{k}= (U^{k}-U^{k-1})/\Delta t, \quad\partial^{2} U^{k}=(\partial U^{k}-\partial U^{k-1})/\Delta t,}\\[2mm]
{\displaystyle \overline{U}^{k} = (U^{k}+U^{k-1})/2, \quad \widetilde{U}^{k}=(U^{k}+U^{k-2})/2,}\\[2mm]
\end{array}
\end{equation}
for any given sequence $\{U^{k}\}_{0}^{M}$ and denote $u^{k}=u(\cdot,t^{k})$ for any given functions
$u\in C(0,T;\,X)$ with a Banach space $ X $.

Let $\mathcal{T}_{h}=\{e\}$ be a regular partition of $\Omega$ into triangles in $\mathbb{R}^{2}$ or tetrahedrons
in $\mathbb{R}^{3}$ without loss of generality, where the mesh size $h= \mathrm{max}_{e\in\mathcal{T}_{h}}\{diam(e)\}$. For a given partition
$\mathcal{T}_{h}$, let $\mathcal{V}_{h}^{r}$ , $\mathbf{V}_{h}^{r}$ and $V_{h}^{r} $ denote the corresponding $r$-th order finite element subspaces of
 $\mathcal{H}_{0}^{1}(\Omega)$ , $\mathbf{H}^{1}_{t}(\Omega)$ and $H_{0}^{1}(\Omega)$, respectively. Let $R_h$, $ {\pi}_h$ and $I_{h}$ be the conventional
 point-wise interpolation operators on $\mathcal{V}_{h}^{r}$, $\mathbf{V}_{h}^{r}$ and $V_{h}^{r}$, respectively.

For convenience, assume that the function $\mathbf{A}$ and $\phi$ is defined in the interval $[-\Delta t, T]$
in terms of the time variable $t$. We can compute $ \mathbf{A}(\cdot,-\Delta t) $ by
\begin{equation}\label{eq:2-4}
{\displaystyle \mathbf{A}(\cdot,-\Delta t)=\mathbf{A}(\cdot,0)-\Delta t \frac{\partial \mathbf{A}}{\partial t}(\cdot,0)=\mathbf{A}_{0}-\Delta t\mathbf{A}_{1}},
\end{equation}
which leads to an approximation to $\mathbf{A}^{-1}$ with second order accuracy. $\phi^{-1}$ can be approximated in a similar way.

An alternating Crank-Nicolson Galerkin finite element approximation to the Maxwell-Schr\"{o}dinger system (\ref{eq:2-1}) is formulated as follows:
\begin{equation}\label{eq:2-5}
\begin{array}{@{}l@{}}
{\displaystyle \psi_{h}^{0}=R_h\psi_{0},\quad\mathbf{A}_{h}^{0}={\pi}_h \mathbf{A}_{0},\quad
\mathbf{A}_{h}^{-1}=\mathbf{A}_{h}^{0}-\Delta t {\pi}_h\mathbf{A}_{1},}\\[2mm]
{\displaystyle \qquad\qquad\phi_{h}^{0} = I_{h}\phi_{0}, \quad \phi_{h}^{-1} = \phi_{h}^{0} - \Delta t I_{h}\psi_{1},}
\end{array}
\end{equation}
and find $(\psi_{h}^{k},\mathbf{A}_{h}^{k},\phi^{k}_{h})\in\mathcal{V}_{h}^{r}\times\mathbf{V}^{r}_{h}\times V^{r}_{h}$ such that for $k=1 ,2 ,\cdots, M$,
\begin{equation}\label{eq:2-6}
\left\{
\begin{array}{@{}l@{}}
{\displaystyle
 -\mathrm{i}({\partial}\psi_{h}^{k},\varphi)+
\frac{1}{2}\left((\mathrm{i}\nabla +\overline{\mathbf{A}}^{k}_{h})\overline{\psi}_{h}^{k},(\mathrm{i}\nabla +\overline{\mathbf{A}}_{h}^{k})\varphi\right)
 +\left( (V_{0}+\overline{\phi}_{h}^{k})\overline{\psi}_{h}^{k},\varphi\right) = 0 ,\,\, \forall \varphi\in\mathcal{V}_{h}^{r},}\\[2mm]
{\displaystyle (\partial^{2}\mathbf{A}_{h}^{k},\mathbf{v}) +\left(\frac{\mathrm{i}}{2}\big((\psi_{h}^{k-1})^{\ast}\nabla{\psi_{h}^{k-1}}
-\psi_h^{k-1}\nabla{(\psi_{h}^{k-1})}^{\ast}\big),\mathbf{v}\right) + (\nabla\times
 \widetilde{\mathbf{A}}_{h}^{k},\nabla\times\mathbf{v})}\\[2mm]
{\displaystyle\quad +(\nabla\cdot
 \widetilde{\mathbf{A}}_{h}^{k},\nabla\cdot\mathbf{v})+\big(|\psi_{h}^{k-1}|^{2}\widetilde{\mathbf{A}}_{h}^{k},\mathbf{v}\big)=0, \quad \forall\mathbf{v}\in\mathbf{V}^{r}_{h},}\\[2mm]
 {\displaystyle (\partial^{2}\phi^{k}_{h}, \eta) + (\nabla \widetilde{\phi}_{h}^{k},\nabla \eta) = (\vert\psi_{h}^{k-1}\vert^{2},\eta),\quad\forall \eta \in V_{h}^{r}.}
\end{array}
\right.
\end{equation}
Note that $ \overline{\mathbf{A}}_{h}^{k} $, $\overline{\phi}_{h}^{k} $,
$ \overline{\psi}_{h}^{k} $, $ \widetilde{\mathbf{A}}_{h}^{k} $ and $ \widetilde{\phi}_{h}^{k} $ are defined
in (\ref{eq:2-2}), and $ (\psi_{h}^{k-1})^{\ast} $ denotes the complex conjugate of $ \psi_{h}^{k-1}$.
For convenience, we define the following bilinear forms:
\begin{equation}\label{eq:2-7}
\begin{array}{@{}l@{}}
{\displaystyle B(\mathbf{A};\psi,\varphi)=\left((\mathrm{i}\nabla+\mathbf{A})\psi,(\mathrm{i}\nabla+\mathbf{A})\varphi\right),}\\[2mm]
{\displaystyle D(\mathbf{A},\mathbf{v})=(\nabla\cdot \mathbf{A},\nabla\cdot \mathbf{v})+
(\nabla\times\mathbf{A},\nabla\times\mathbf{v}),}\\[2mm]
{\displaystyle f(\psi,\varphi)=\frac{\mathrm{i}}{2}(\varphi^{\ast}\nabla\psi-\psi\nabla\varphi^{\ast}).}
\end{array}
\end{equation}

Then the variational form of the Maxwell-Schr\"{o}dinger system (\ref{eq:2-1}) and its discrete system (\ref{eq:2-2})  can be reformulated as follows:
\begin{equation}\label{eq:2-8}
\left\{
\begin{array}{@{}l@{}}
{\displaystyle -\mathrm{i}(\frac{\partial \psi}{\partial t},\varphi)+\frac12B(\mathbf{A};\psi,\varphi)+(V_{0}\psi,\varphi) + (\phi\psi,\varphi)=0,\quad \forall \varphi\in\mathcal{H}_{0}^{1}(\Omega),}\\[2mm]
{\displaystyle (\frac{\partial ^{2}\mathbf{A}}{\partial t^{2}},\mathbf{v})
+D(\mathbf{A},\mathbf{v})+( f(\psi,\psi),\mathbf{v})+( |\psi|^{2}\mathbf{A},\mathbf{v})=0, \quad \forall\mathbf{v}\in\mathbf{H}^{1}_{t}(\Omega),}\\[2mm]
{\displaystyle  (\frac{\partial^{2}\phi}{\partial t^{2}},\eta)+(\nabla \phi,\nabla\eta)=(\vert\psi\vert^{2},\eta), \quad\quad \forall \eta\in H_{0}^{1}(\Omega)}
\end{array}
\right.
\end{equation}
and for $ k=1,2, \cdots,M$,
\begin{equation}\label{eq:2-9}
\left\{
\begin{array}{@{}l@{}}
{\displaystyle
 -\mathrm{i}({\partial}\psi_{h}^{k},\varphi)+
\frac{1}{2}B(\overline{\mathbf{A}}_{h}^{k};\overline{\psi}_{h}^{k},\varphi)
 +( V_{0}\overline{\psi}_{h}^{k},\varphi)+(\overline{\phi}_{h}^{k}\overline{\psi}_{h}^{k},\varphi)=0 ,\quad \forall \varphi\in\mathcal{V}_{h}^{r},}\\[2mm]
{\displaystyle (\partial ^{2}\mathbf{A}_{h}^{k},\mathbf{v})+D(\widetilde{\mathbf{A}}_{h}^{k},\mathbf{v})
+\left( f(\psi_{h}^{k-1},\psi_h^{k-1}),\mathbf{v}\right)+\big( |\psi_{h}^{k-1}|^{2}\widetilde{\mathbf{A}}_{h}^{k},\mathbf{v}\big) = 0, \quad
 \forall\mathbf{v}\in\mathbf{V}^{r}_{h},}\\[2mm]
  {\displaystyle (\partial^{2}\phi^{k}_{h}, \eta) + (\nabla \widetilde{\phi}_{h}^{k},\nabla \eta) = (\vert\psi_{h}^{k-1}\vert^{2},\eta),\quad\forall \eta \in V_{h}^{r}.}
\end{array}
\right.
\end{equation}

In this paper we assume that the Maxwell-Schr\"{o}dinger equations (\ref{eq:2-8}) has one and only one weak solution $(\psi,\mathbf{A},\phi)$ and
the following regularity conditions are satisfied:
\begin{equation}\label{eq:2-10}
\begin{array}{@{}l@{}}
{\displaystyle
\psi,\psi_{t},\psi_{tt} \in {L}^{\infty}(0, T; \mathcal{H}^{r+1}(\Omega)),\quad \psi_{ttt} \in {L}^{\infty}(0, T; \mathcal{H}^{1}(\Omega)),}\\[2mm]
{\displaystyle\psi_{tttt} \in L^{2}(0, T; \mathcal{L}^{2}(\Omega));\quad
\mathbf{A},\mathbf{A}_{t},\mathbf{A}_{tt} \in {L}^{\infty}(0, T; \mathbf{H}^{r+1}(\Omega)),}\\[2mm]
{\displaystyle\mathbf{A}_{ttt} \in {L}^{\infty}(0, T; \mathbf{H}^{1}(\Omega)),\mathbf{A}_{tttt} \in L^{2}(0, T; \mathbf{L}^{2}(\Omega));}\\[2mm]
{\displaystyle \phi,\phi_{t},\phi_{tt}\in {L}^{\infty}(0, T; {H}^{r+1}(\Omega)), \phi_{ttt}\in {L}^{\infty}(0, T; {H}^{1}(\Omega)),}\\[2mm]
{\displaystyle \phi_{tttt}\in L^{2}(0, T; {L}^{2}(\Omega)).}
\end{array}
\end{equation}

For the initial conditions $(\psi_{0},\mathbf{A}_{0},\mathbf{A}_{1},\phi_{0},\phi_{1})$ , we assume that
\begin{equation}\label{eq:2-11}
\begin{array}{@{}l@{}}
{\displaystyle \psi_{0}\in \mathcal{H}^{r+1}(\Omega) \cap \mathcal{H}_{0}^{1}(\Omega); \,\,\, \mathbf{A}_{0},\mathbf{A}_{1}\in\mathbf{H}^{r+1}(\Omega)\cap\mathbf{H}_{t}^{1}(\Omega);\,\,
\, \phi_{0},\phi_{1} \in H^{r+1}(\Omega)\cap H_{0}^{1}(\Omega).}
\end{array}
\end{equation}

We now give the main convergence result in this paper as follows:
\begin{theorem}\label{thm2-1}
Suppose that the Maxwell-Schr\"{o}dinger coupled system (\ref{eq:2-8}) has a unique solution
$(\psi,\mathbf{A},\phi)$ satisfying (\ref{eq:2-10}) and (\ref{eq:2-11}).
Let $(\psi_h^k,\mathbf{A}_h^k,\phi_h^k)$ be the fully discrete numerical solution
of $(\psi,\mathbf{A},\phi)$ defined in (\ref{eq:2-9}). 
 Then there exist two positive constants $h_{0}>0$
and $\Delta t_{0}>0$, such that when $h < h_{0}$, $\Delta t < \Delta t_{0}$, we have the following error estimates:
\begin{equation}\label{eq:2-12}
\begin{array}{@{}l@{}}
 {\displaystyle \max_{1\leq k \leq M}\big\{\|\psi_{h}^{k}-\psi^{k}\|_{\mathcal{L}^2(\Omega)}^{2}
 +\|\nabla(\psi_{h}^{k}-\psi^{k})\|_{\mathbf{L}^2(\Omega)}^{2}
 + \|\mathbf{A}_{h}^{k}-\mathbf{A}^{k}\|_{\mathbf{L}^2(\Omega)}^{2}}\\[2mm]
{\displaystyle \qquad\quad+\|\nabla\cdot(\mathbf{A}_{h}^{k}-\mathbf{A}^{k})\|_{{L}^2(\Omega)}^{2}
+\|\nabla\times(\mathbf{A}_{h}^{k}-\mathbf{A}^{k})\|_{\mathbf{L}^2(\Omega)}^{2} }\\[2mm]
{\displaystyle \qquad + \|\phi_{h}^{k}-\phi^{k}\|^{2}_{L^{2}(\Omega)}
+\|\nabla(\phi_{h}^{k}-\phi^{k})\|^{2}_{\mathbf{L}^{2}(\Omega)} \big\}\leq C \big\{h^{2r}+(\Delta t)^{4}\big\},\,\, r\geq 1,}
\end{array}
\end{equation}
where $ \psi^k=\psi(\cdot, t^k) $, $ \mathbf{A}^k=\mathbf{A}(\cdot, t^k) $, $ \phi^k=\phi(\cdot, t^k) $, and $C$ is a constant independent of $h$, $\Delta t$.
\end{theorem}

\section{The proof of Theorem~\ref{thm2-1}}\label{sec-3}
In this section, we will give the proof of Theorem~\ref{thm2-1}.

\subsection{Preliminaries}
For convenience, we list some imbedding inequalities and interpolation inequalities in Sobolev spaces (see, e.g., \cite{Lad} and \cite{Gir}),
and use them in the sequel.
\begin{equation}\label{eq:3-1}
 {\displaystyle \|u\|_{L^p} \leq C \|u\|_{{H}^1}, \quad \|\mathbf{v}\|_{\mathbf{L}^p}\leq C \|\mathbf{v}\|_{\mathbf{H}^1},
 \quad 1\leq p \leq 6 \,\,(d=2,3),}
 \end{equation}
 \begin{equation}\label{eq:3-2}
 {\displaystyle  \|\mathbf{v}\|_{\mathbf{H}^1} \leq C(\|\nabla\times\mathbf{v}\|_{\mathbf{L}^2}+\|\nabla\cdot\mathbf{v}\|_{{L}^2}),\quad
 \mathbf{v}\in \mathbf{H}^{1}_{t}(\Omega),}
 \end{equation}
 \begin{equation}\label{eq:3-2-1}
 {\displaystyle \|u\|_{{L}^3}^{2}\leq \|u\|_{{L}^2}\|u\|_{{L}^6},}
 \end{equation}
where $ \|u\|_{L^p}=\|u\|_{L^p(\Omega)} $, $ \|u\|_{{H}^1}=\|u\|_{{H}^1(\Omega)} $,
$ \|\mathbf{v}\|_{\mathbf{L}^2}=\|\mathbf{v}\|_{\mathbf{L}^2(\Omega)} $ and
$\|\mathbf{v}\|_{\mathbf{H}^1}=\|\mathbf{v}\|_{\mathbf{H}^1(\Omega)} $.

The following identities will be used frequently in this paper.
\begin{equation}\label{eq:3-2-2}
\begin{array}{@{}l@{}}
{\displaystyle \sum_{k=1}^{M}{(a_{k}-a_{k-1})b_{k}}=a_{M}b_{M}-a_{0}b_{1}-\sum_{k=1}^{M-1}{a_{k}(b_{k+1}-b_{k})},}\\[2mm]
{\displaystyle \sum_{k=1}^{M}{(a_{k}-a_{k-1})b_{k}}=a_{M}b_{M}-a_{0}b_{0}-\sum_{k=1}^{M}{a_{k-1}(b_{k}-b_{k-1})}.}
\end{array}
\end{equation}

Let $(R_h\psi, {\pi}_{h}\mathbf{A}, I_{h}\phi)$ denote the interpolation function of $(\psi, \mathbf{A},\phi)$
in $\mathcal{V}_{h}^{r}\times\mathbf{V}_{h}^{r}\times V_{h}^{r}$. Set
$e_{\psi}=R_h\psi-\psi$,  $e_{\mathbf{A}}={\pi}_h\mathbf{A}-\mathbf{A}$, $e_{\phi} = I_{h}\phi-\phi $. By applying standard finite element theory and the regualrity conditions in (\ref{eq:2-10}), we have
\begin{equation}\label{eq:3-2-3}
\begin{array}{@{}l@{}}
{\displaystyle \|e_{\psi}\|_{\mathcal{L}^2}+h\|e_{\psi}\|_{\mathcal{H}^1}\leq Ch^{r+1},\quad  \|e_{\mathbf{A}}\|_{\mathbf{L}^2}
+h\|e_{\mathbf{A}}\|_{\mathbf{H}^1}\leq C h^{r+1},}\\[2mm]
{\displaystyle \|e_{\phi}\|_{{L}^2}+h\|e_{\phi}\|_{{H}^1}\leq Ch^{r+1},\,\,\,
\|R_h\psi\|_{\mathcal{L}^{\infty}}+\|{\pi}_h\mathbf{A}\|_{\mathbf{L}^{\infty}}+\|I_{h}\phi\|_{{L}^{\infty}}
+\|\nabla R_h\psi \|_{\mathbf{L}^{3}}\leq C,}
\end{array}
\end{equation}
where $ C $ is a constant independent of $ h $.

The following lemmas will be useful in the proof of Theorem~\ref{thm2-1}.
\begin{lemma}\label{lem3-1}
For the solution of (\ref{eq:2-9}), we have
\begin{equation}\label{eq:3-2-4}
{\displaystyle \|\psi_{h}^{k}\|_{\mathcal{L}^2}^{2}={\|\psi_{h}^{0}\|}_{\mathcal{L}^2}^{2},\quad k=1,2\cdots, M.}
\end{equation}
\end{lemma}

\begin{lemma}\label{lem3-2}
For $k=1,2\cdots,M $, the following identities hold for the bilinear functional $B(\mathbf{A};\psi,\varphi)$ defined in (\ref{eq:2-7}):
\begin{equation}\label{eq:3-2-5}
\begin{array}{@{}l@{}}
{\displaystyle B(\mathbf{A};\psi,\varphi)-B(\widetilde{\mathbf{A}};\psi,\varphi)=\left((\mathbf{A}
+\widetilde{\mathbf{A}})\psi \varphi^{*},\mathbf{A}-\widetilde{\mathbf{A}}\right)
+2(f(\psi,\varphi),\mathbf{A}-\widetilde{\mathbf{A}}),}\\[2mm]
{\displaystyle \mathrm{Re}\left[B\left(\overline{\mathbf{A}}^{k};\overline{\psi}^{k},\partial {\psi}^{k}\right)\right]=
 -\left(\frac{1}{2}(\overline{\mathbf{A}}^{k}+\overline{\mathbf{A}}^{k-1})|{\psi}^{k-1}|^{2},\frac{1}{2}(\partial \mathbf{A}^{k}+\partial \mathbf{A}^{k-1})\right)}\\[2mm]
{\displaystyle\quad -\left(f({\psi}^{k-1}, {\psi}^{k-1}),\frac{1}{2}(\partial \mathbf{A}^{k}+\partial \mathbf{A}^{k-1})\right) +\frac{1}{2}\partial B(\overline{\mathbf{A}}^{k};{\psi}^{k},{\psi}^{k}).}
\end{array}
\end{equation}
\end{lemma}

Lemma~\ref{lem3-1} can be proved by choosing $\varphi=\overline{\psi}_{h}^{k}$ in $(\ref{eq:2-9})_1 $ and taking the imaginary part. A direct calculation gives (\ref{eq:3-2-5}) in Lemma~\ref{lem3-2}.

Let $\theta_{\psi}^{k}=\psi_{h}^{k}-R_h\psi^{k} $, $\theta_{\mathbf{A}}^{k}=\mathbf{A}_{h}^{k}-{\pi}_{h}\mathbf{A}^{k}$, $\theta_{\phi}^{k}= \phi_{h}^{k}-I_{h}\phi^{k}$.
By using the error estimates of the interpolation operators (\ref{eq:3-2-3}), we only need to estimate $\theta_{\psi}^{k}$ , $\theta_{\mathbf{A}}^{k}$ and $\theta_{\phi}^{k}$. We will prove the following estimate:
\begin{equation}\label{eq:3-2-6}
{\displaystyle\|\theta_{\psi}^{k}\|^{2}_{\mathcal{H}^{1}} + \|\partial \theta_{\mathbf{A}}^{k}\|^{2}_{\mathbf{L}^{2}} +  \|\theta_{\mathbf{A}}^{k}\|^{2}_{\mathbf{H}^{1}} + \|\partial\theta_{\phi}^{k}\|^{2}_{L^2} +
\|\theta_{\phi}^{k}\|^{2}_{H^1} \leq C_{\ast}\left\{h^{2r}+(\Delta t)^4\right\},}
\end{equation}
where $k=0,1,\cdots,M$.

By using the regularity assumption of the initial conditions (\ref{eq:2-11}) and the error estimates of the interpolation operators (\ref{eq:3-2-3}), we get
\begin{equation}\label{eq:3-2-7}
{\displaystyle\|\theta_{\psi}^{0}\|^{2}_{\mathcal{H}^{1}} + \|\partial \theta_{\mathbf{A}}^{0}\|^{2}_{\mathbf{L}^{2}} +  \|\theta_{\mathbf{A}}^{0}\|^{2}_{\mathbf{H}^{1}} + \|\partial\theta_{\phi}^{0}\|^{2}_{L^2} +
\|\theta_{\phi}^{0}\|^{2}_{H^1} \leq C_{0}\left\{h^{2r}+(\Delta t)^4\right\}.}
\end{equation}

We use the mathematical inductive method to show (\ref{eq:3-2-6}). By (\ref{eq:3-2-7}), if we require $C_{\ast} \geq C_{0} $, then (\ref{eq:3-2-6}) holds for $k=0$. We assume that (\ref{eq:3-2-6}) holds for $0\leq k \leq m-1$. In the rest of this section,
we will find $C_{\ast}$ such that (\ref{eq:3-2-6}) holds for $k\leq m$, where $C_{\ast}$ is independent of $k\,$, $h\,$, $\Delta t$.

Subtracting (\ref{eq:2-8}) from (\ref{eq:2-9}), we obtain the following equations for $\theta_{\mathbf{A}}^{k}$, $\theta_{\phi}^{k}$ and $\theta_{\psi}^{k}$ :
\begin{equation}\label{eq:3-3}
\begin{array}{@{}l@{}}
{\displaystyle \left(\partial^{2}\theta^{k}_{\mathbf{A}},\mathbf{v}\right)+D(\widetilde{\theta}^{k}_{\mathbf{A}},\mathbf{v})
= \left(\frac{\partial^{2}\mathbf{A}^{k-1}}{\partial t^{2}}-\partial^{2} {\pi}_{h}\mathbf{A}^{k},\mathbf{v}\right)}\\[2mm]
{\displaystyle \quad +D(\mathbf{A}^{k-1}-\widetilde{{\pi}_{h}\mathbf{A}}^{k},\mathbf{v})
+ \left(|\psi^{k-1}|^{2}\mathbf{A}^{k-1}-|\psi^{k-1}_{h}|^{2}\widetilde{\mathbf{A}}^{k}_{h},
\;\mathbf{v}\right)}\\[2mm]
{\displaystyle \quad+ \left(f(\psi^{k-1},\psi^{k-1})-f(\psi^{k-1}_{h},\psi^{k-1}_{h}),\;\mathbf{v}\right),\quad \forall\mathbf{v}\in\mathbf{V}^{r}_{h},}
\end{array}
\end{equation}
\begin{equation}\label{eq:3-3-0}
\begin{array}{@{}l@{}}
{\displaystyle \left(\partial^{2}\theta^{k}_{\phi},\eta\right)+\left(\nabla \widetilde{\theta}^{k}_{\phi},\nabla \eta\right) = \left(\frac{\partial^{2}\phi^{k-1}}{\partial t^{2}}-\partial^{2}I_{h}\phi^{k},\eta\right)}\\[2mm]
{\displaystyle + \left(\nabla(\phi^{k-1}-\widetilde{I_{h}\phi}^{k}),\nabla \eta\right) + \left(\vert\psi^{k-1}_{h}\vert^{2}-\vert\psi^{k-1}\vert^{2},\eta\right),\quad \forall \eta \in V_{h}^{r},}
\end{array}
\end{equation}
\begin{equation}\label{eq:3-4}
\begin{array}{@{}l@{}}
{\displaystyle-2\mathrm{i}(\partial \theta^{k}_{\psi},\varphi)+B\left(\overline{\mathbf{A}}^{k}_{h};\overline{\theta}_{\psi}^{k},\varphi\right)
= 2\mathrm{i}\left(\partial R_h\psi^{k}-\frac{\partial \psi^{k-\frac{1}{2}}}{\partial t},\varphi\right)+2V_0\left(\psi^{k-\frac{1}{2}}-\overline{\psi}_{h}^{k},\varphi\right)}\\[2mm]
{\displaystyle \quad+ B(\mathbf{A}^{k-\frac{1}{2}};(\psi^{k-\frac{1}{2}}-R_h\overline{\psi}^{k}),\varphi)
+2\left(\phi^{k-\frac12}\psi^{k-\frac12}-\overline{\phi}^{k}_{h}\overline{\psi}^{k}_{h},
\varphi\right)}\\[2mm]
{\displaystyle \quad+\left(B(\mathbf{A}^{k-\frac{1}{2}};R_h\overline{\psi}^{k},\varphi)-B(\overline{\mathbf{A}}^{k}_{h}; R_h\overline{\psi}^{k},\varphi)\right),\quad \forall \varphi\in\mathcal{V}_{h}^{r}.}
\end{array}
\end{equation}

The key steps of the proof of (\ref{eq:3-2-6}) are now briefly described. In order to find $C_{\ast} $ and to show that (\ref{eq:3-2-6}) holds for $k=m$, we first take $\mathbf{v}
=\frac{\displaystyle 1}{\displaystyle 2\Delta t}(\theta_{\mathbf{A}}^{k}-\theta_{\mathbf{A}}^{k-2})$ in (\ref{eq:3-3}) and $\eta = \frac{\displaystyle 1}{\displaystyle 2\Delta t}(\theta_{\phi}^{k}-\theta_{\phi}^{k-2}) $ in (\ref{eq:3-3-0}) and obtain the estimates of $\theta_{\mathbf{A}}^m $ and $\theta_{\phi}^m $. Then we choose $\varphi=\overline{\theta}_{\psi}^{k}$ in (\ref{eq:3-4}) and give the estimate of $ \|\theta_\psi^m\|_{\mathcal{L}^2} $. Finally, we take $\varphi=\partial \theta_{\psi}^{k}$ in (\ref{eq:3-4})
 and make use of the above estimates of $\theta_{\mathbf{A}}^m $ and $\theta_{\phi}^m $ to derive the energy-norm estimate for $\theta^{m}_{\psi}$.
 Using the above estimates, we can complete the proof of (\ref{eq:3-2-6}).

\subsection{Estimates for (\ref{eq:3-3})}
If we set
\begin{equation}\label{eq:3-61}
\begin{array}{@{}l@{}}
{\displaystyle  I_1^{k}(\mathbf{v})=\left(\frac{\partial^{2}\mathbf{A}^{k-1}}{\partial t^{2}}-\partial^{2}{\pi}_{h}\mathbf{A}^{k},\mathbf{v}\right),\quad I_2^{k}(\mathbf{v})= D(\mathbf{A}^{k-1}-\widetilde{{\pi}_{h}\mathbf{A}}^{k},\mathbf{v}),}\\[2mm]
{\displaystyle I_3^{k}(\mathbf{v}) = \left(|\psi^{k-1}|^{2}\mathbf{A}^{k-1}-|\psi^{k-1}_{h}|^{2}\widetilde{\mathbf{A}}^{k}_{h},\;\mathbf{v}\right),}\\[2mm]
{\displaystyle I_4^{k}(\mathbf{v})= \left(f(\psi^{k-1},\psi^{k-1})-f(\psi^{k-1}_{h},\psi^{k-1}_{h}),\;\mathbf{v}\right),\quad \mathbf{v}\in\mathbf{V}_{h}^{r},}
\end{array}
\end{equation}
then we rewrite (\ref{eq:3-3}) as follows:
\begin{equation}\label{eq:3-62}
{\displaystyle \left(\partial^{2}\theta^{k}_{\mathbf{A}},\mathbf{v}\right)+D(\widetilde{\theta}^{k}_{\mathbf{A}},\mathbf{v})=
 I_1^{k}(\mathbf{v})+ I_2^{k}(\mathbf{v})+ I_3^{k}(\mathbf{v})+ I_4^{k}(\mathbf{v}).}
\end{equation}

We take $\mathbf{v}= \frac{\displaystyle 1}{\displaystyle 2\Delta t}(\theta_{\mathbf{A}}^{k}
-\theta_{\mathbf{A}}^{k-2})=\frac{\displaystyle 1}{\displaystyle 2}(\partial \theta_{\mathbf{A}}^{k}
+\partial \theta_{\mathbf{A}}^{k-1}) = \overline{\partial \theta}_{\mathbf{A}}^{k}$ in (\ref{eq:3-62}) and get
\begin{equation}\label{eq:3-62-0}
\begin{array}{@{}l@{}}
{\displaystyle \left(\partial^{2}\theta^{k}_{\mathbf{A}},\frac{1}{2}(\partial \theta_{\mathbf{A}}^{k}+\partial \theta_{\mathbf{A}}^{k-1})\right)
+D(\widetilde{\theta^{k}_{\mathbf{A}}},\frac{1}{2}(\partial \theta_{\mathbf{A}}^{k}+\partial \theta_{\mathbf{A}}^{k-1}))}\\[2mm]
{\displaystyle \quad =\frac{1}{2\varDelta t}\left[\|\partial \theta_{\mathbf{A}}^{k}\|_{\mathbf{L}^2}^{2}
-\|\partial \theta_{\mathbf{A}}^{k-1}\|_{\mathbf{L}^2}^{2}\right]+\frac{1}{4\Delta t}\left[D({\theta}_{\mathbf{A}}^{k}, {\theta}_{\mathbf{A}}^{k})
-D({\theta}_{\mathbf{A}}^{k-2}, {\theta}_{\mathbf{A}}^{k-2})\right]}\\[2mm]
{\displaystyle  \quad =I_1^{k}(\overline{\partial \theta}_{\mathbf{A}}^{k})+ I_2^{k}(\overline{\partial \theta}_{\mathbf{A}}^{k})+ I_3^{k}(\overline{\partial \theta}_{\mathbf{A}}^{k}) + I_4^{k}(\overline{\partial \theta}_{\mathbf{A}}^{k}),}
\end{array}
\end{equation}
which leads to
\begin{equation}\label{eq:3-62-1}
\begin{array}{@{}l@{}}
{\displaystyle  \frac{1}{2}\|\partial \theta_{\mathbf{A}}^{m}\|_{\mathbf{L}^2}^{2}
+\frac{1}{4}D({\theta}_{\mathbf{A}}^{m}, {\theta}_{\mathbf{A}}^{m})
+\frac{1}{4}D({\theta}_{\mathbf{A}}^{m-1}, {\theta}_{\mathbf{A}}^{m-1})}\\[2mm]
{\displaystyle \quad= \frac{1}{2}\|\partial \theta_{\mathbf{A}}^{0}\|_{\mathbf{L}^2}^{2}
+ \frac{1}{4}D({\theta}_{\mathbf{A}}^{0}, {\theta}_{\mathbf{A}}^{0})
+\frac{1}{4}D({\theta}_{\mathbf{A}}^{-1}, {\theta}_{\mathbf{A}}^{-1})}\\[2mm]
{\displaystyle\quad+{\Delta t}\sum_{k=1}^{m} \left[I_1^{k}(\overline{\partial \theta}_{\mathbf{A}}^{k})+ I_2^{k}(\overline{\partial \theta}_{\mathbf{A}}^{k})+ I_3^{k}(\overline{\partial \theta}_{\mathbf{A}}^{k})+ I_4^{k}(\overline{\partial \theta}_{\mathbf{A}}^{k})\right]}\\[2mm]
{\displaystyle \quad \leq Ch^{2r} +{\Delta t}\sum_{k=1}^{m} \left[I_1^{k}(\overline{\partial \theta}_{\mathbf{A}}^{k})+ I_2^{k}(\overline{\partial \theta}_{\mathbf{A}}^{k})+ I_3^{k}(\overline{\partial \theta}_{\mathbf{A}}^{k}) + I_4^{k}(\overline{\partial \theta}_{\mathbf{A}}^{k})\right].}
\end{array}
\end{equation}
Here we have used the fact that
\begin{equation*}
{\displaystyle \|\partial \theta_{\mathbf{A}}^{0}\|_{\mathbf{L}^2}^{2}
+D({\theta}_{\mathbf{A}}^{0}, {\theta}_{\mathbf{A}}^{0})
+D({\theta}_{\mathbf{A}}^{-1}, {\theta}_{\mathbf{A}}^{-1})
\leq C h^{2r}.}
\end{equation*}

Now we estimate  $\sum\limits_{k=1}^{m} I_j^{k}(\overline{\partial \theta}_{\mathbf{A}}^{k})$, $j=1,2,3,4$.
Under the regularity assumption of $\mathbf{A}$ in (\ref{eq:2-10}), we can prove
\begin{equation}\label{eq:3-63}
\begin{array}{@{}l@{}}
{\displaystyle\sum_{k=1}^{m} |I_1^{k}(\overline{\partial \theta}_{\mathbf{A}}^{k})| \leq \sum_{k=1}^{m}|\left[\frac{\partial^{2}\mathbf{A}^{k-1}}{\partial t^{2}}-\partial^{2}\mathbf{A}^{k},\;\overline{\partial \theta}_{\mathbf{A}}^{k}\right]|}\\[2mm]
{\displaystyle \quad +\sum_{k=1}^{m}
|\left(\partial^{2}\mathbf{A}^{k}-{\pi}_{h}\partial^{2}\mathbf{A}^{k},\;\overline{\partial \theta}_{\mathbf{A}}^{k}\right)|}\\[2mm]
{\displaystyle \quad\leq \frac{C}{\Delta t}\left\{h^{2r+2}+(\Delta t)^{4}\right\}
+ C\sum_{k=1}^{m} \|\overline{\partial \theta}_{\mathbf{A}}^{k}\|_{\mathbf{L}^2}^{2}.}
\end{array}
\end{equation}

We rewrite the term $\Delta t \sum_{k=1}^{m} I_2^{k}(\overline{\partial \theta}_{\mathbf{A}}^{k})$ as follows:
\begin{equation}\label{eq:3-63-1}
\begin{array}{@{}l@{}}
{\displaystyle  \Delta t \sum_{k=1}^{m} I_2^{k}(\overline{\partial \theta}_{\mathbf{A}}^{k})= \frac{1}{2}\Delta t \sum_{k=1}^{m} D\left(\mathbf{A}^{k-1}
-\frac{1}{2}(\mathbf{A}^{k}+\mathbf{A}^{k-2}), (\partial \theta_{\mathbf{A}}^{k}
+\partial \theta_{\mathbf{A}}^{k-1})\right)}\\[2mm]
{\displaystyle \quad + \frac{1}{4}\Delta t \sum_{k=1}^{m} D\left((
\mathbf{A}^{k}+\mathbf{A}^{k-2})-{\pi}_{h}(
\mathbf{A}^{k}+\mathbf{A}^{k-2}), (\partial \theta_{\mathbf{A}}^{k}+\partial \theta_{\mathbf{A}}^{k-1})\right).}
\end{array}
\end{equation}

Applying (\ref{eq:3-2-2}), the regularity assumption and Young's inequality, we get
\begin{equation}\label{eq:3-63-2}
\begin{array}{@{}l@{}}
{\displaystyle \frac{1}{2}\Delta t\sum_{k=1}^{m} D\left(\mathbf{A}^{k-1}
-\frac{1}{2}(\mathbf{A}^{k}+\mathbf{A}^{k-2}), (\partial \theta_{\mathbf{A}}^{k}
+\partial \theta_{\mathbf{A}}^{k-1})\right)}\\[2mm]
{\displaystyle \leq C\left\{h^{2r} + (\Delta t)^{4}\right\}
+\frac{1}{64} D\left(\theta_{\mathbf{A}}^{m},\theta_{\mathbf{A}}^{m}\right)
+\frac{1}{64} D\left(\theta_{\mathbf{A}}^{m-1},\theta_{\mathbf{A}}^{m-1}\right) + C\Delta t \sum_{k=0}^{m}D\left(\theta_{\mathbf{A}}^{k},\theta_{\mathbf{A}}^{k}\right)}\\[2mm]
\end{array}
\end{equation}
and
\begin{equation}\label{eq:3-63-3}
\begin{array}{@{}l@{}}
{\displaystyle \frac{1}{4}\Delta t \sum_{k=1}^{m}D\left((
\mathbf{A}^{k}+\mathbf{A}^{k-2})-{\pi}_{h}(
\mathbf{A}^{k}+\mathbf{A}^{k-2}), (\partial \theta_{\mathbf{A}}^{k}+\partial \theta_{\mathbf{A}}^{k-1})\right)}\\[2mm]
{\displaystyle \leq C\left\{h^{2r} + (\Delta t)^{4}\right\}
+\frac{1}{64} D\left(\theta_{\mathbf{A}}^{m},\theta_{\mathbf{A}}^{m}\right)
+\frac{1}{64} D\left(\theta_{\mathbf{A}}^{m-1},\theta_{\mathbf{A}}^{m-1}\right)
+C\Delta t \sum_{k=0}^{m}D\left(\theta_{\mathbf{A}}^{k},\theta_{\mathbf{A}}^{k}\right).}
\end{array}
\end{equation}

We thus have
\begin{equation}\label{eq:3-63-4}
\begin{array}{@{}l@{}}
{\displaystyle  \Delta t \sum_{k=1}^{m} I_2^{k}(\overline{\partial \theta}_{\mathbf{A}}^{k})\leq C\left\{h^{2r} + (\Delta t)^{4}\right\}
+ C\Delta t \sum_{k=0}^{m}D\left(\theta_{\mathbf{A}}^{k},\theta_{\mathbf{A}}^{k}\right).}\\[2mm]
{\displaystyle \quad\quad+ \frac{1}{32} D\left(\theta_{\mathbf{A}}^{m},\theta_{\mathbf{A}}^{m}\right)
+\frac{1}{32} D\left(\theta_{\mathbf{A}}^{m-1},\theta_{\mathbf{A}}^{m-1}\right).}
\end{array}
\end{equation}

We observe that
\begin{equation}\label{eq:3-64}
\begin{array}{@{}l@{}}
{\displaystyle I_3^{k}(\overline{\partial \theta}_{\mathbf{A}}^{k}) = \left(|\psi^{k-1}|^{2}(\mathbf{A}^{k-1}-{\pi}_{h}\widetilde{\mathbf{A}}^{k}),\;\overline{\partial \theta}_{\mathbf{A}}^{k}\right)+\left(|\psi^{k-1}|^{2}({\pi}_{h}\widetilde{\mathbf{A}}^{k}-\widetilde{\mathbf{A}}^{k}_{h}),\;\overline{\partial \theta}_{\mathbf{A}}^{k}\right) }\\[2mm]
{\displaystyle \quad + \left((|\psi^{k-1}|^{2} - |R_h\psi^{k-1}|^{2}){\pi}_{h}\widetilde{\mathbf{A}}^{k},\;\overline{\partial \theta}_{\mathbf{A}}^{k}\right) + \left((|\psi^{k-1}|^{2} - |R_h\psi^{k-1}|^{2})\widetilde{\theta}_\mathbf{A}^{k},\;\overline{\partial \theta}_{\mathbf{A}}^{k}\right) }\\[2mm]
{\displaystyle \quad\ +\left((|R_h\psi^{k-1}|^{2} - |\psi^{k-1}_{h}|^{2}){\pi}_{h}\widetilde{\mathbf{A}}^{k},\;\overline{\partial \theta}_{\mathbf{A}}^{k}\right) +\left((|R_h\psi^{k-1}|^{2} - |\psi^{k-1}_{h}|^{2})\widetilde{\theta}_\mathbf{A}^{k},\;\overline{\partial \theta}_{\mathbf{A}}^{k}\right)  }\\[2mm]
{\displaystyle \quad\stackrel{\mathrm{def}}{=}\sum_{j=1}^{6} I_3^{k,j}(\overline{\partial \theta}_{\mathbf{A}}^{k}) .}
\end{array}
\end{equation}

The first four terms in (\ref{eq:3-64}) can be estimated by a standard argument, i.e.
\begin{equation}\label{eq:3-65}
\begin{array}{@{}l@{}}
{\displaystyle \sum_{j=1}^{4}|I_3^{k,j}(\overline{\partial \theta}_{\mathbf{A}}^{k})| \leq C \left\{(\Delta t)^{4}+ h^{2r+2}+D(\widetilde{\theta}_{\mathbf{A}}^{k}, \widetilde{\theta}_{\mathbf{A}}^{k})+\|\overline{\partial \theta}_{\mathbf{A}}^{k}\|_{\mathbf{L}^2}^{2}\right\}.}
\end{array}
\end{equation}

We notice that
\begin{equation}\label{eq:3-65-0}
\begin{array}{@{}l@{}}
{\displaystyle |R_h\psi^{k-1}|^{2} - |\psi^{k-1}_{h}|^{2} = (R_h\psi^{k-1}-\psi^{k-1}_{h})(R_h\psi^{k-1})^{\ast} + (R_h\psi^{k-1}-\psi^{k-1}_{h})^{\ast} \psi^{k-1}_{h} } \\[2mm]
{\displaystyle \quad\quad\qquad\qquad\qquad\quad= -\theta_{\psi}^{k-1}(R_h\psi^{k-1})^{\ast}
-(\theta_{\psi}^{k-1})^{\ast}R_h\psi^{k-1}-\vert\theta_{\psi}^{k-1}\vert^{2}}
\end{array}
\end{equation}
and obtain
 \begin{equation}\label{eq:3-66}
 \begin{array}{@{}l@{}}
 {\displaystyle \vert I_3^{k,5}(\overline{\partial \theta}_{\mathbf{A}}^{k})\vert  \leq C\|\theta_{\psi}^{k-1}\|_{\mathcal{L}^6}  \|\overline{\partial \theta}_{\mathbf{A}}^{k}\|_{\mathbf{L}^{2}}  +  C \|\theta_{\psi}^{k-1}\|^{2}_{\mathcal{L}^6} \|\overline{\partial \theta}_{\mathbf{A}}^{k}\|_{\mathbf{L}^{2}},  }\\[2mm]
 {\displaystyle \vert I_3^{k,6}(\overline{\partial \theta}_{\mathbf{A}}^{k})\vert  \leq C\|\theta_{\psi}^{k-1}\|_{\mathcal{L}^6}  \| \widetilde{\theta}_{\mathbf{A}}^{k} \|_{\mathbf{L}^{6}} \|\overline{\partial \theta}_{\mathbf{A}}^{k}\|_{\mathbf{L}^{2}}  + C \|\theta_{\psi}^{k-1}\|^{2}_{\mathcal{L}^6}\|\widetilde{\theta}_{\mathbf{A}}^{k} \|_{\mathbf{L}^{6}} \|\overline{\partial \theta}_{\mathbf{A}}^{k}\|_{\mathbf{L}^{2}}.}
 \end{array}
 \end{equation}

 By using the assumption of the induction, we have
 \begin{equation*}
  {\displaystyle\|\theta_{\psi}^{k-1}\|_{\mathcal{L}^6} \leq C \|\theta_{\psi}^{k-1}\|_{\mathcal{H}^1} \leq CC^{\frac12}_{\ast}\left\{(\Delta t)^{2}+h ^{r}\right\}.}
 \end{equation*}

 If we choose some sufficiently small $ h $ and
 $ \Delta t $ such that $ CC^{\frac12}_{\ast}\left\{(\Delta t)^{2}+h ^{r}\right\} \leq 1$, then we get
 \begin{equation}\label{eq:3-67}
  \begin{array}{@{}l@{}}
  {\displaystyle \vert I_3^{k,5}(\overline{\partial \theta}_{\mathbf{A}}^{k})\vert  \leq C \|\theta_{\psi}^{k-1}\|_{\mathcal{L}^6}\|\overline{\partial \theta}_{\mathbf{A}}^{k}\|_{\mathbf{L}^{2}} \leq C\left\{ \|\nabla \theta_{\psi}^{k-1}\|^{2}_{\mathbf{L}^2} + \|\overline{\partial \theta}_{\mathbf{A}}^{k}\|^{2}_{\mathbf{L}^{2}}\right\}, } \\[2mm]
 {\displaystyle \vert I_3^{k,6}(\overline{\partial \theta}_{\mathbf{A}}^{k})\vert  \leq C \|\widetilde{\theta}_{\mathbf{A}}^{k}\|_{\mathbf{L}^6}\|\overline{\partial \theta}_{\mathbf{A}}^{k}\|_{\mathbf{L}^{2}}\leq C \left\{D(\widetilde{\theta}_{\mathbf{A}}^{k}, \widetilde{\theta}_{\mathbf{A}}^{k})+\|\overline{\partial \theta}_{\mathbf{A}}^{k}\|^{2}_{\mathbf{L}^{2}} \right\}.}
 \end{array}
 \end{equation}

Combining (\ref{eq:3-64})-(\ref{eq:3-67}) implies
 \begin{equation}\label{eq:3-68}
 \begin{array}{@{}l@{}}
 {\displaystyle \sum_{k=1}^{m} \vert I_3^{k}(\overline{\partial \theta}_{\mathbf{A}}^{k}) \vert \leq \frac{C}{\Delta t}\left\{h^{2r}+(\Delta t)^{4}\right\}+ C\sum_{k=1}^{m}\left\{D(\widetilde{\theta}_{\mathbf{A}}^{k}, \widetilde{\theta}_{\mathbf{A}}^{k})+\|\nabla\theta_{\psi}^{k-1}\|_{\mathbf{L}^2}^{2}
 +\|\overline{\partial \theta}_{\mathbf{A}}^{k}\|_{\mathbf{L}^2}^{2}\right\}.}
 \end{array}
\end{equation}

To estimate  $\sum\limits_{k=1}^{m} I_4^{k}(\mathbf{v})$, we rewrite
$ I_4^{k}(\overline{\partial \theta}_{\mathbf{A}}^{k})$ as follows:
\begin{equation}\label{eq:3-69}
\begin{array}{@{}l@{}}
{\displaystyle I_4^{k}(\overline{\partial \theta}_{\mathbf{A}}^{k})=\left(f(\psi^{k-1},\psi^{k-1})-f(R_h\psi^{k-1}, R_h\psi^{k-1}),\;\overline{\partial \theta}_{\mathbf{A}}^{k}\right) }\\[2mm]
{\displaystyle \quad +\left(f(R_h\psi^{k-1}, R_h\psi^{k-1})-f(\psi^{k-1}_{h},\psi^{k-1}_{h}),\;\overline{\partial \theta}_{\mathbf{A}}^{k}\right)\stackrel{\mathrm{def}}{=} I_4^{k,1}(\overline{\partial \theta}_{\mathbf{A}}^{k})+I_4^{k,2}(\overline{\partial \theta}_{\mathbf{A}}^{k}).}
\end{array}
\end{equation}

Observing
\begin{equation}\label{eq:3-70}
\begin{array}{@{}l@{}}
{\displaystyle f(\psi,\psi)-f(\varphi,\varphi)=\frac{\mathrm{i}}{2}(\psi^{\ast}\nabla \psi- \psi \nabla\psi^{\ast})
  -\frac{\mathrm{i}}{2} (\varphi^{\ast} \nabla\varphi-\varphi \nabla\varphi^{\ast}) }\\[2mm]
{\displaystyle = -\frac{\mathrm{i}}{2}\left(\varphi^{\ast}\nabla(\varphi-\psi)-\varphi\nabla(\varphi-\psi)^{\ast}\right)
 +\frac{\mathrm{i}}{2}\left((\varphi-\psi)\nabla\psi^{\ast}-(\varphi-\psi)^{\ast}\nabla\psi\right)}
\end{array}
\end{equation}
and applying (\ref{eq:3-2-3}), we get
\begin{equation}\label{eq:3-71}
\begin{array}{@{}l@{}}
{\displaystyle I_4^{k,1}(\overline{\partial \theta}_{\mathbf{A}}^{k})=\left(f(\psi^{k-1},\psi^{k-1})-f(R_h\psi^{k-1}, R_h\psi^{k-1}),\;\overline{\partial \theta}_{\mathbf{A}}^{k}\right)}\\[2mm]
{\displaystyle \quad = -\frac{\mathrm{i}}{2}\left((R_h\psi^{k-1})^{\ast}\nabla(R_h\psi^{k-1}-\psi^{k-1})
-R_h\psi^{k-1}\nabla(R_h\psi^{k-1}-\psi^{k-1})^{\ast},\;\overline{\partial \theta}_{\mathbf{A}}^{k}\right)}\\[2mm]
{\displaystyle \quad +\frac{\mathrm{i}}{2}\left((R_h\psi^{k-1}-\psi^{k-1})\nabla(\psi^{k-1})^{\ast}-(R_h\psi^{k-1}-\psi^{k-1})^{\ast}
\nabla\psi^{k-1},\;\overline{\partial \theta}_{\mathbf{A}}^{k}\right)}\\[2mm]
{\displaystyle\quad \leq  C\left\{h^{2r}+\|\overline{\partial \theta}_{\mathbf{A}}^{k}\|_{\mathbf{L}^2}^{2}\right\}.}
\end{array}
\end{equation}

Using (\ref{eq:3-2-3}), we similarly prove
\begin{equation}\label{eq:3-72}
\begin{array}{@{}l@{}}
{\displaystyle I_4^{k,2}(\overline{\partial \theta}_{\mathbf{A}}^{k})=\left(f(R_h\psi^{k-1}, R_h\psi^{k-1})-f(\psi^{k-1}_{h},\psi^{k-1}_{h}),\;\overline{\partial \theta}_{\mathbf{A}}^{k}\right)}\\[2mm]
{\displaystyle \quad = -\frac{\mathrm{i}}{2}\left((\psi^{k-1}_{h})^{\ast}\nabla(\psi^{k-1}_{h}-R_h\psi^{k-1})
-\psi^{k-1}_{h}\nabla(\psi^{k-1}_{h}-R_h\psi^{k-1})^{\ast},\;\overline{\partial \theta}_{\mathbf{A}}^{k}\right)}\\[2mm]
{\displaystyle \quad\quad +\frac{\mathrm{i}}{2}\left((\psi^{k-1}_{h}-R_h\psi^{k-1})\nabla (R_h\psi^{k-1})^{\ast}-(\psi^{k-1}_{h}-R_h\psi^{k-1})^{\ast}\nabla R_h\psi^{k-1},\;\overline{\partial \theta}_{\mathbf{A}}^{k}\right)}\\[2mm]
{\displaystyle \quad = -\frac{\mathrm{i}}{2}\left((\theta_{\psi}^{k-1})^{\ast}\nabla\theta_{\psi}^{k-1}
-\theta_{\psi}^{k-1}\nabla(\theta_{\psi}^{k-1})^{\ast},\;\overline{\partial \theta}_{\mathbf{A}}^{k}\right) }\\[2mm]
{\displaystyle \quad\quad-\frac{\mathrm{i}}{2}\left((R_h\psi^{k-1})^{\ast}\nabla\theta_{\psi}^{k-1}
-R_h\psi^{k-1}\nabla(\theta_{\psi}^{k-1})^{\ast},\;\overline{\partial \theta}_{\mathbf{A}}^{k}\right) }\\[2mm]
{\displaystyle\quad \quad +\frac{\mathrm{i}}{2}\left(\theta_{\psi}^{k-1}\nabla (R_h\psi^{k-1})^{\ast}-(\theta_{\psi}^{k-1})^{\ast}\nabla R_h\psi^{k-1},\;\overline{\partial \theta}_{\mathbf{A}}^{k}\right)}\\[2mm]
{\displaystyle \quad \leq -\left(f(\theta_{\psi}^{k-1},\theta_{\psi}^{k-1}),\;\overline{\partial \theta}_{\mathbf{A}}^{k}\right)
+ C \| R_h\psi^{k-1}\|_{\mathcal{L}^{\infty}}\|\nabla\theta_{\psi}^{k-1}\|_{\mathbf{L}^2}\|\overline{\partial \theta}_{\mathbf{A}}^{k}\|_{\mathbf{L}^2}}\\[2mm]
{\displaystyle \quad\quad+C\|\nabla R_h\psi^{k-1}\|_{\mathbf{L}^{3}}\|\theta_{\psi}^{k-1}\|_{\mathcal{L}^6}
\|\overline{\partial \theta}_{\mathbf{A}}^{k}\|_{\mathbf{L}^2}}\\[2mm]
{\displaystyle \quad \leq -\left(f(\theta_{\psi}^{k-1},\theta_{\psi}^{k-1}),\;\overline{\partial \theta}_{\mathbf{A}}^{k}\right)
+C\left\{\|\nabla\theta_{\psi}^{k-1}\|_{\mathbf{L}^2}^{2}
+\|\overline{\partial \theta}_{\mathbf{A}}^{k}\|_{\mathbf{L}^2}^{2}\right\}.}
\end{array}
\end{equation}

Hence we have
\begin{equation}\label{eq:3-73}
\begin{array}{@{}l@{}}
{\displaystyle\sum_{k=1}^{m} I_4^{k}(\overline{\partial \theta}_{\mathbf{A}}^{k})\leq \frac{C h^{2r}}{\Delta t} -\sum_{k=1}^{m}\left(f(\theta_{\psi}^{k-1},\theta_{\psi}^{k-1}),\;\overline{\partial \theta}_{\mathbf{A}}^{k}\right)
+C\sum_{k=1}^{m}\left\{\|\nabla\theta_{\psi}^{k-1}\|_{\mathbf{L}^2}^{2}
+\|\overline{\partial \theta}_{\mathbf{A}}^{k}\|_{\mathbf{L}^2}^{2}\right\}.}
\end{array}
\end{equation}

Substituting (\ref{eq:3-63}), (\ref{eq:3-63-4}), (\ref{eq:3-68}) and (\ref{eq:3-73}) into (\ref{eq:3-62-1}),  we get
\begin{equation}\label{eq:3-76}
\begin{array}{@{}l@{}}
{\displaystyle  \frac{1}{2}\|\partial \theta_{\mathbf{A}}^{m}\|_{\mathbf{L}^2}^{2}
+\frac{7}{32}D({\theta}_{\mathbf{A}}^{m}, {\theta}_{\mathbf{A}}^{m})
+\frac{7}{32}D({\theta}_{\mathbf{A}}^{m-1}, {\theta}_{\mathbf{A}}^{m-1})}\\[2mm]
{\displaystyle \quad \leq C\left\{h^{2r}+(\Delta t)^{4} \right\}
+ C\Delta t\sum_{k=0}^{m}\left(D({\theta}_{\mathbf{A}}^{k}, {\theta}_{\mathbf{A}}^{k})
+ \|\partial \theta_{\mathbf{A}}^{k}\|_{\mathbf{L}^2}^{2} \right)}\\[2mm]
{\displaystyle \quad +C\Delta t\sum_{k=0}^{m-1} \|\nabla\theta_{\psi}^{k}\|_{\mathbf{L}^2}^{2}-\Delta t \sum_{k=1}^{m}\left(f(\theta_{\psi}^{k-1},\theta_{\psi}^{k-1}),\;\overline{\partial \theta}_{\mathbf{A}}^{k}\right).}
\end{array}
\end{equation}

We estimate the last term on the right side of (\ref{eq:3-76}). It follows from the definition of bilinear functional $f(\varphi, \psi)$ in (\ref{eq:2-7}) that
\begin{equation}
\begin{array}{@{}l@{}}
{\displaystyle \sum_{k=1}^{m}\left(f(\theta_{\psi}^{k-1},\theta_{\psi}^{k-1}),\overline{\partial \theta}_{\mathbf{A}}^{k}\right) = \frac{\mathrm{i}}{4}\sum_{k=1}^{m}\left((\theta_{\psi}^{k-1})^{\ast}\nabla\theta_{\psi}^{k-1}
-\theta_{\psi}^{k-1}\nabla(\theta_{\psi}^{k-1})^{\ast},\; (\partial \theta_{\mathbf{A}}^{k}
+\partial \theta_{\mathbf{A}}^{k-1})\right).}
\end{array}
\end{equation}

For $\Delta t \leq h^{\frac12}$, by the assumption of the induction and the inverse inequalities,  we have
\begin{equation*}
{\displaystyle\|\nabla\theta_{\psi}^{k-1}\|_{\mathbf{L}^{3}} \leq C h^{-\frac12}\|\nabla\theta_{\psi}^{k-1}\|_{\mathbf{L}^{2}} \leq C h^{-\frac12}C_{\ast}^{\frac12}\{h^{r}+h\}\leq CC_{\ast}^{\frac12}h^{\frac12}.}
\end{equation*}

We choose a sufficiently small $ h>0 $ such that $CC_{\ast}^{\frac12}h^{\frac12}\leq 1 $ and  obtain
\begin{equation*}
{\displaystyle\|\nabla\theta_{\psi}^{k-1}\|_{\mathbf{L}^{3}} \leq 1.}
\end{equation*}
Consequently
\begin{equation}\label{eq:3-81}
\begin{array}{@{}l@{}}
{\displaystyle \vert \sum_{k=1}^{m}\left(f(\theta_{\psi}^{k-1},\theta_{\psi}^{k-1}),\overline{\partial \theta}_{\mathbf{A}}^{k}\right)\vert \leq  \frac 12 \sum_{k=1}^{m} \| \theta_{\psi}^{k-1} \|_{\mathcal{L}^{6}} \|\nabla\theta_{\psi}^{k-1}\|_{\mathbf{L}^{3}} \|\partial \theta_{\mathbf{A}}^{k}
+\partial \theta_{\mathbf{A}}^{k-1}\|_{\mathbf{L}^{2}}    }\\[2mm]
{\displaystyle \qquad \leq  \frac 12 \sum_{k=1}^{m} \| \theta_{\psi}^{k-1} \|_{\mathcal{L}^{6}}\|\partial \theta_{\mathbf{A}}^{k}
+\partial \theta_{\mathbf{A}}^{k-1}\|_{\mathbf{L}^{2}}  \leq C \sum_{k=1}^{m} \| \nabla \theta_{\psi}^{k-1} \|_{\mathbf{L}^{2}}\|\partial \theta_{\mathbf{A}}^{k}
+\partial \theta_{\mathbf{A}}^{k-1}\|_{\mathbf{L}^{2}}  }\\[2mm]
{\displaystyle \qquad \leq C \sum_{k=0}^{m-1} \| \nabla \theta_{\psi}^{k} \|^{2}_{\mathbf{L}^{2}}
+ C \sum_{k=0}^{m} \|\partial \theta_{\mathbf{A}}^{k} \|^{2}_{\mathbf{L}^{2}},}
\end{array}
\end{equation}
where $ C $ is a constant independent of $ h $, $ \Delta t $.

As for $h^{\frac12} \leq \Delta t$, by the assumption of the induction, we discover
\begin{equation*}
{\displaystyle\|\nabla\theta_{\psi}^{k-1}\|_{\mathbf{L}^{2}} \leq C_{\ast}^{\frac12}\left\{(\Delta t)^{2} + (\Delta t)^{2r}\right\} \leq 2C_{\ast}^{\frac12}(\Delta t)^{2}.}
\end{equation*}

Now choose a sufficiently small $ \Delta t>0 $ to find $2C_{\ast}^{\frac12}\Delta t \leq 1$, in which case we have
\begin{equation*}
{\displaystyle\|\nabla\theta_{\psi}^{k-1}\|_{\mathbf{L}^{2}} \leq \Delta t.}
\end{equation*}
It follows that
\begin{equation}\label{eq:3-82}
\begin{array}{@{}l@{}}
{\displaystyle \vert \sum_{k=1}^{m}\left(f(\theta_{\psi}^{k-1},\theta_{\psi}^{k-1}),\overline{\partial \theta}_{\mathbf{A}}^{k}\right) \vert \leq  \frac 12 \sum_{k=1}^{m} \| \theta_{\psi}^{k-1} \|_{\mathcal{L}^{3}} \|\nabla\theta_{\psi}^{k-1}\|_{\mathbf{L}^{2}} \|\partial \theta_{\mathbf{A}}^{k}
+\partial \theta_{\mathbf{A}}^{k-1}\|_{\mathbf{L}^{6}} }\\[2mm]
{\displaystyle \qquad \leq \frac {\Delta t}{2} \sum_{k=1}^{m} \| \theta_{\psi}^{k-1} \|_{\mathcal{L}^{3}} \|\partial \theta_{\mathbf{A}}^{k}
+\partial \theta_{\mathbf{A}}^{k-1}\|_{\mathbf{L}^{6}}  \leq \frac {1}{2} \sum_{k=1}^{m} \| \theta_{\psi}^{k-1} \|_{\mathcal{L}^{3}} \|\theta_{\mathbf{A}}^{k}
-\theta_{\mathbf{A}}^{k-2}\|_{\mathbf{L}^{6}} }\\[2mm]
{\displaystyle \qquad \leq C \sum_{k=1}^{m} \| \nabla \theta_{\psi}^{k-1} \|^{2}_{\mathbf{L}^{2}}
+ C\sum_{k=1}^{m} (D(\theta_{\mathbf{A}}^{k},\theta_{\mathbf{A}}^{k}) + D(\theta_{\mathbf{A}}^{k-2},\theta_{\mathbf{A}}^{k-2}))  }\\[2mm]
{\displaystyle \qquad \leq C h^{2r} + C \sum_{k=0}^{m-1} \| \nabla \theta_{\psi}^{k} \|^{2}_{\mathbf{L}^{2}}
+ C\sum_{k=0}^{m} D(\theta_{\mathbf{A}}^{k},\theta_{\mathbf{A}}^{k}). }
\end{array}
\end{equation}

Combining (\ref{eq:3-81}), (\ref{eq:3-82}) and (\ref{eq:3-76}) implies
\begin{equation}\label{eq:3-83}
 \begin{array}{@{}l@{}}
{\displaystyle  \frac{1}{2}\|\partial \theta_{\mathbf{A}}^{m}\|_{\mathbf{L}^2}^{2}
+\frac{7}{32}D({\theta}_{\mathbf{A}}^{m}, {\theta}_{\mathbf{A}}^{m})
+\frac{7}{32}D({\theta}_{\mathbf{A}}^{m-1}, {\theta}_{\mathbf{A}}^{m-1}) \leq C\left\{h^{2r}+(\Delta t)^{4} \right\} }\\[2mm]
 {\displaystyle \qquad \quad + C\Delta t\sum_{k=0}^{m}\left(D({\theta}_{\mathbf{A}}^{k}, {\theta}_{\mathbf{A}}^{k}) + \|\partial \theta_{\mathbf{A}}^{k}\|_{\mathbf{L}^2}^{2} \right) +C\Delta t\sum_{k=0}^{m-1} \|\nabla\theta_{\psi}^{k}\|_{\mathbf{L}^2}^{2}. }
 \end{array}
 \end{equation}

 We choose a sufficiently small $ \Delta t>0 $ such that $C\Delta t \leq \frac{1}{8}$ and find
 \begin{equation}
  \begin{array}{@{}l@{}}
{\displaystyle  \frac{3}{8}\|\partial \theta_{\mathbf{A}}^{m}\|_{\mathbf{L}^2}^{2}
+\frac{3}{32}D({\theta}_{\mathbf{A}}^{m}, {\theta}_{\mathbf{A}}^{m})
+\frac{3}{32}D({\theta}_{\mathbf{A}}^{m-1}, {\theta}_{\mathbf{A}}^{m-1}) \leq C\left\{h^{2r}+(\Delta t)^{4} \right\} }\\[2mm]
 {\displaystyle \qquad \quad + C\Delta t\sum_{k=0}^{m-1}\left(D({\theta}_{\mathbf{A}}^{k}, {\theta}_{\mathbf{A}}^{k}) + \|\partial \theta_{\mathbf{A}}^{k}\|_{\mathbf{L}^2}^{2} \right) +C\Delta t\sum_{k=0}^{m-1} \|\nabla\theta_{\psi}^{k}\|_{\mathbf{L}^2}^{2}, }
 \end{array}
 \end{equation}
 which leads to
  \begin{equation}\label{eq:3-84}
  \begin{array}{@{}l@{}}
{\displaystyle  \frac{3}{8}\|\partial \theta_{\mathbf{A}}^{m}\|_{\mathbf{L}^2}^{2}
+\frac{3}{32}D({\theta}_{\mathbf{A}}^{m}, {\theta}_{\mathbf{A}}^{m})
+\frac{3}{32}D({\theta}_{\mathbf{A}}^{m-1}, {\theta}_{\mathbf{A}}^{m-1})  }\\[2mm]
 {\displaystyle \qquad \quad  \leq C\left\{h^{2r}+(\Delta t)^{4} \right\} +  CC_{\ast}\left\{h^{2r}+(\Delta t)^{4} \right\} }
 \end{array}
 \end{equation}
by the assumption of the induction.

Applying the assumption of the induction again and (\ref{eq:3-84}), for $k=1,2,\cdots,m$, we deduce
\begin{equation}\label{eq:3-84-0}
  \begin{array}{@{}l@{}}
  {\displaystyle \|\partial \mathbf{A}_{h}^{k} \|_{\mathbf{L}^{2}} + \|\mathbf{A}_{h}^{k}\|_{\mathbf{L}^{6}} \leq \|\partial \theta_{\mathbf{A}}^{k}\|_{\mathbf{L}^2} + \|\theta_{\mathbf{A}}^{k}\|_{\mathbf{L}^{6}}
  +\|\partial{\pi}_{h} \mathbf{A}^{k} \|_{\mathbf{L}^{2}}  + \|{\pi}_{h}\mathbf{A}^{k}\|_{\mathbf{L}^{6}} } \\[2mm]
  {\displaystyle \qquad\quad \leq C\left\{\|\partial \theta_{\mathbf{A}}^{k}\|_{\mathbf{L}^2}+ D^{\frac12}({\theta}_{\mathbf{A}}^{k}, {\theta}_{\mathbf{A}}^{k})\right\}+\|\partial{\pi}_{h} \mathbf{A}^{k} \|_{\mathbf{L}^{2}}+
  \|{\pi}_{h}\mathbf{A}^{k}\|_{\mathbf{L}^{6}}\leq C, }
\end{array}
\end{equation}
when $\|\partial \theta_{\mathbf{A}}^{k}\|_{\mathbf{L}^2}+ D^{\frac12}({\theta}_{\mathbf{A}}^{k}, {\theta}_{\mathbf{A}}^{k}) \leq 1.$

\subsection{Estimates for (\ref{eq:3-3-0})}
 Setting
 \begin{equation}\label{eq:3-85}
 \begin{array}{@{}l@{}}
{\displaystyle J_1^{k}(\eta) =  \left(\frac{\partial^{2}\phi^{k-1}}{\partial t^{2}}-\partial^{2}I_{h}\phi^{k},\,\eta\right) ,}\\[2mm]
{\displaystyle J_2^{k}(\eta) = \left(\nabla(\phi^{k-1}-\widetilde{I_{h}\phi}^{k}),\,\nabla \eta\right), }\\[2mm]
 {\displaystyle J_3^{k}(\eta) = \left(\vert\psi^{k-1}_{h}\vert^{2}-\vert\psi^{k-1}\vert^{2},\,\eta\right),\quad  \eta \in V_{h}^{r},}
\end{array}
 \end{equation}
 we rewrite (\ref{eq:3-3-0}) as follows:
 \begin{equation}\label{eq:3-86}
  \begin{array}{@{}l@{}}
{\displaystyle  \left(\partial^{2}\theta^{k}_{\phi},\,\eta\right)+\left(\nabla \widetilde{\theta}^{k}_{\phi},\,\nabla \eta\right)  =J_1^{k}(\eta) + J_2^{k}(\eta) + J_3^{k}(\eta).}
 \end{array}
 \end{equation}

 We take $\eta = \frac{\displaystyle 1}{\displaystyle 2}(\partial \theta_{\phi}^{k}+\partial \theta_{\phi}^{k-1})
 =\overline{\partial \theta}_{\phi}^{k}$ in (\ref{eq:3-86}) and obtain
 \begin{equation}\label{eq:3-87}
  \begin{array}{@{}l@{}}
 {\displaystyle  \left(\partial^{2}\theta^{k}_{\phi},\,\overline{\partial \theta}_{\phi}^{k}\right)+\left(\nabla \widetilde{\theta}^{k}_{\phi},\,\nabla \overline{\partial \theta}_{\phi}^{k}\right) }\\[2mm]
 {\displaystyle\,\,  = \frac{1}{2\Delta t} \left(\|\partial \theta^{k}_{\phi}\|^{2}_{{L}^{2}} - \|\partial \theta^{k-1}_{\phi}\|^{2}_{{L}^{2}}\right) + \frac{1}{4\Delta t} \left(\|\nabla \theta^{k}_{\phi}\|^{2}_{\mathbf{L}^{2}} - \|\nabla \theta^{k-2}_{\phi}\|^{2}_{\mathbf{L}^{2}}\right) }\\[2mm]
 {\displaystyle \,\, =J_1^{k}(\overline{\partial \theta}_{\phi}^{k}) + J_2^{k}(\overline{\partial \theta}_{\phi}^{k}) + J_3^{k}(\overline{\partial \theta}_{\phi}^{k}).}
 \end{array}
 \end{equation}

 Multiply (\ref{eq:3-87}) by $\Delta t$ and sum $k =1, 2,\cdots,m  $ to discover
 \begin{equation}\label{eq:3-88}
 \begin{array}{@{}l@{}}
 {\displaystyle \frac{1}{2}\|\partial \theta^{m}_{\phi}\|^{2}_{{L}^{2}} + \frac{1}{4}\|\nabla \theta^{m}_{\phi}\|^{2}_{\mathbf{L}^{2}} + \frac{1}{4}\|\nabla \theta^{m-1}_{\phi}\|^{2}_{\mathbf{L}^{2}} } \\[2mm]
 {\displaystyle\quad = \frac{1}{2}\|\partial \theta^{0}_{\phi}\|^{2}_{{L}^{2}} + \frac{1}{4}\|\nabla \theta^{0}_{\phi}\|^{2}_{\mathbf{L}^{2}} + \frac{1}{4}\|\nabla \theta^{-1}_{\phi}\|^{2}_{\mathbf{L}^{2}}}\\[2mm]
{\displaystyle \quad + \Delta t \sum_{k=1}^{m} \left\{J_1^{k}(\overline{\partial \theta}_{\phi}^{k}) + J_2^{k}(\overline{\partial \theta}_{\phi}^{k}) + J_3^{k}(\overline{\partial \theta}_{\phi}^{k})\right\}} \\[2mm]
{\displaystyle \quad \leq  C\big\{h^{2r} + (\Delta t)^{4} \big\} +  \Delta t \sum_{k=1}^{m} \left\{J_1^{k}(\overline{\partial \theta}_{\phi}^{k}) + J_2^{k}(\overline{\partial \theta}_{\phi}^{k}) + J_3^{k}(\overline{\partial \theta}_{\phi}^{k})\right\}.}
 \end{array}
 \end{equation}

 Applying the regularity assumption of $\phi $ in (\ref{eq:2-10}), we deduce
 \begin{equation}\label{eq:3-89}
  \begin{array}{@{}l@{}}
  {\displaystyle\sum_{k=1}^{m} |J_1^{k}(\overline{\partial \theta}_{\phi}^{k})| \leq \sum_{k=1}^{m}|\left(\frac{\partial^{2}{\phi}^{k-1}}{\partial t^{2}}-\partial^{2}{\phi}^{k},\,\overline{\partial \theta}_{\phi}^{k}\right)|+\sum_{k=1}^{m}
|\left(\partial^{2}{\phi}^{k}-I_{h}\partial^{2}{\phi}^{k},\,\overline{\partial \theta}_{\phi}^{k}\right)|}\\[2mm]
{\displaystyle \quad\quad \quad\leq \frac{C}{\Delta t}\big\{h^{2r+2}+(\Delta t)^{4}\big\}
+ C\sum_{k=1}^{m} \|\overline{\partial \theta}_{\phi}^{k}\|_{{L}^2}^{2}.}
 \end{array}
 \end{equation}
 $\sum_{k=1}^{m} J_2^{k}(\overline{\partial \theta}_{\phi}^{k})$ can be bounded by a argument similar to the estimate of  $\sum_{k=1}^{m} I_2^{k}(\overline{\partial \theta}_{\mathbf{A}}^{k})$
 in (\ref{eq:3-63-2}) and (\ref{eq:3-63-3}).
 \begin{equation}\label{eq:3-90}
 \begin{array}{@{}l@{}}
{\displaystyle  \Delta t \sum_{k=1}^{m} J_2^{k}(\overline{\partial \theta}_{{\phi}}^{k})\leq C\left\{h^{2r} + (\Delta t)^{4}\right\}
+ \frac{1}{8} \|\nabla \theta^{m}_{\phi}\|^{2}_{\mathbf{L}^{2}}}\\[2mm]
{\displaystyle\quad+\frac{1}{8}\|\nabla \theta^{m-1}_{\phi}\|^{2}_{\mathbf{L}^{2}} + C\Delta t \sum_{k=0}^{m}\|\nabla \theta^{k}_{\phi}\|^{2}_{\mathbf{L}^{2}}. }
\end{array}
\end{equation}

Observing
\begin{equation*}
{\displaystyle J_3^{k}(\overline{\partial \theta}_{\phi}^{k})=\left(\vert R_h\psi^{k-1}\vert^{2}-\vert\psi^{k-1}\vert^{2},\, \overline{\partial \theta}_{\phi}^{k}\right)+ \left(\vert \psi^{k-1}_{h}\vert^{2}-\vert R_h\psi^{k-1}\vert^{2},\, \overline{\partial \theta}_{\phi}^{k}\right),}
\end{equation*}
and using (\ref{eq:3-65-0}), we get
 \begin{equation}\label{eq:3-91}
 \begin{array}{@{}l@{}}
{\displaystyle  \Delta t \sum_{k=1}^{m} J_3^{k}(\overline{\partial \theta}_{{\phi}}^{k})\leq C h^{2r} + C\Delta t \sum_{k=1}^{m} \left\{\|\nabla \theta_{\psi}^{k-1}\|^{2}_{\mathbf{L}^{2}} + \|\overline{\partial \theta}_{\phi}^{k}\|_{{L}^2}^{2} \right\}}\\[2mm]
{\displaystyle \quad + \Delta t \sum_{k=1}^{m}\|\theta_{\psi}^{k-1}\|^{2}_{\mathcal{L}^{4}}\|\overline{\partial \theta}_{\phi}^{k}\|_{{L}^2}
\leq C h^{2r} + C\Delta t \sum_{k=1}^{m} \left\{\|\nabla \theta_{\psi}^{k-1}\|^{2}_{\mathbf{L}^{2}} + \|\overline{\partial \theta}_{\phi}^{k}\|_{{L}^2}^{2} \right\}} \\[2mm]
{\displaystyle \quad + \Delta t \sum_{k=1}^{m}\|\theta_{\psi}^{k-1}\|_{\mathcal{L}^{4}}\|\overline{\partial \theta}_{\phi}^{k}\|_{{L}^2}\leq C h^{2r} + C\Delta t \sum_{k=1}^{m} \left\{\|\nabla \theta_{\psi}^{k-1}\|^{2}_{\mathbf{L}^{2}} + \|\overline{\partial \theta}_{\phi}^{k}\|_{{L}^2}^{2} \right\}.}
 \end{array}
 \end{equation}
 Here we have used $\|\theta_{\psi}^{k-1}\|_{\mathcal{L}^{4}}  \leq  C \|\theta_{\psi}^{k-1}\|_{\mathcal{H}^{1}} \leq 1.$

 Substituting (\ref{eq:3-89}), (\ref{eq:3-90}) and (\ref{eq:3-91}) into (\ref{eq:3-88}) yields
  \begin{equation}\label{eq:3-92}
 \begin{array}{@{}l@{}}
 {\displaystyle \frac{1}{2}\|\partial \theta^{m}_{\phi}\|^{2}_{{L}^{2}} + \frac{1}{8}\|\nabla \theta^{m}_{\phi}\|^{2}_{\mathbf{L}^{2}} + \frac{1}{8}\|\nabla \theta^{m-1}_{\phi}\|^{2}_{\mathbf{L}^{2}} \leq C\left\{h^{2r} + (\Delta t)^{4}\right\} } \\[2mm]
 {\displaystyle  \qquad + C\Delta t \sum_{k=0}^{m} \left(\|\nabla \theta_{\phi}^{k}\|^{2}_{\mathbf{L}^{2}} + \|{\partial \theta}_{\phi}^{k}\|_{{L}^2}^{2} \right) + C\Delta t \sum_{k=0}^{m-1} \|\nabla \theta_{\psi}^{k}\|^{2}_{\mathbf{L}^{2}}. }
 \end{array}
 \end{equation}

 Similarly to (\ref{eq:3-84-0}), for $k=1,2,\cdots,m$, we can prove
 \begin{equation}\label{eq:3-93}
{\displaystyle
 \|\partial \phi_{h}^{k}\|_{L^{2}} + \|\phi_{h}^{k}\|_{{L}^{6}} \leq C,}
 \end{equation}
 where $ C $ is a constant independent of $ h $, $ \Delta t $ and $C_{\ast}$.

\subsection{Estimates for (\ref{eq:3-4})}

We rewrite (\ref{eq:3-4}) as follows:
 \begin{equation}\label{eq:3-5}
 {\displaystyle -2\mathrm{i}(\partial \theta^{k}_{\psi},\varphi)+B\left(\overline{\mathbf{A}}^{k}_{h};\overline{\theta}_{\psi}^{k},\varphi\right)
 =\sum\limits_{j=1}^5 Q_j^{k}(\varphi),}
 \end{equation}
 where
 \begin{equation*}
 \begin{array}{@{}l@{}}
 {\displaystyle Q_1^{k}(\varphi)=2\mathrm{i}\left(\partial R_h\psi^{k}-\frac{\partial \psi^{k-\frac{1}{2}}}{\partial t},\,\varphi\right),\quad Q_2^{k}(\varphi)=2V_0\left(\psi^{k-\frac{1}{2}}-\overline{\psi}_{h}^{k},\,\varphi\right),}\\[2mm]
 {\displaystyle Q_3^{k}(\varphi)=B(\mathbf{A}^{k-\frac{1}{2}};(\psi^{k-\frac{1}{2}}-R_h\overline{\psi}^{k}),\varphi), \quad Q_4^{k}(\varphi) = 2\left(\phi^{k-\frac12}\psi^{k-\frac12}-\overline{\phi}^{k}_{h}\overline{\psi}^{k}_{h},\,\varphi\right), } \\[2mm]
 {\displaystyle Q_5^{k}(\varphi)=B(\mathbf{A}^{k-\frac{1}{2}};R_h\overline{\psi}^{k},\varphi)-B(\overline{\mathbf{A}}^{k}_{h};
 R_h\overline{\psi}^{k},\varphi).}
 \end{array}
 \end{equation*}

 Taking $\varphi =\overline{\theta}_{\psi}^{k}$ in (\ref{eq:3-5}), and observing the imaginary part of the above equation, we have
 \begin{equation}\label{eq:3-5-0}
  \begin{array}{@{}l@{}}
 {\displaystyle \frac{1}{\Delta t} \left(\| {\theta}_{\psi}^{k} \|^{2}_{\mathcal{L}^{2}} - \| {\theta}_{\psi}^{k-1} \|^{2}_{\mathcal{L}^{2}}\right) = -\mathrm{Im}\Big[Q_1^{k}( \overline{\theta}_{\psi}^{k})+Q_2^{k}( \overline{\theta}_{\psi}^{k})+Q_3^{k}( \overline{\theta}_{\psi}^{k})+Q_4^{k}( \overline{\theta}_{\psi}^{k})}\\[2mm]
{\displaystyle \quad + Q_5^{k}( \overline{\theta}_{\psi}^{k}) \Big]\leq \vert Q_1^{k}( \overline{\theta}_{\psi}^{k}) \vert +\vert Q_2^{k}( \overline{\theta}_{\psi}^{k}) \vert + \vert Q_3^{k}( \overline{\theta}_{\psi}^{k}) \vert + \vert Q_4^{k}( \overline{\theta}_{\psi}^{k}) \vert
  +\vert Q_5^{k}( \overline{\theta}_{\psi}^{k}) \vert.}
 \end{array}
 \end{equation}

 It is obvious that
\begin{equation*}
 \begin{array}{@{}l@{}}
{\displaystyle  Q_1^{k}( \overline{\theta}_{\psi}^{k})=2\mathrm{i}\left(\partial \psi^{k}-\frac{\partial \psi^{k-\frac{1}{2}}}{\partial t},\,\overline{\theta}_{\psi}^{k}\right)+\left(R_h\partial\psi^{k}-\partial \psi^{k},\,\overline{\theta}_{\psi}^{k}\right),}\\[2mm]
{\displaystyle  Q_2^{k}( \overline{\theta}_{\psi}^{k})= -2V_0\left(\overline{\theta}_{\psi}^{k},\overline{\theta}_{\psi}^{k}\right)+2V_0\left(  \overline{\psi}^{k} - R_h\overline{\psi}^{k}, \,\overline{\theta}_{\psi}^{k}\right)+ 2V_0\left(\psi^{k-\frac{1}{2}}- \overline{\psi}^{k},\,\overline{\theta}_{\psi}^{k}\right).}
\end{array}
\end{equation*}

Using the error estimates (\ref{eq:3-2-3}) for the interpolation operator $R_h$ and
the regularity of $\psi$ in (\ref{eq:2-10}), we give the following estimate
 \begin{equation}\label{eq:3-6}
 \begin{array}{@{}l@{}}
{\displaystyle  |Q_1^{k}( \overline{\theta}_{\psi}^{k})|+ |Q_2^{k}( \overline{\theta}_{\psi}^{k})|\leq C\big\{(\Delta t)^{4}+h^{2r+2}\big\}+C\|\overline{\theta}_{\psi}^{k}\|^{2}_{\mathcal{L}^2}.}
\end{array}
\end{equation}

Note that
\begin{equation}\label{eq:3-7}
\begin{array}{@{}l@{}}
 {\displaystyle B(\mathbf{A};\psi,\varphi)=(\nabla\psi,\nabla\varphi)+\left(|\mathbf{A}|^{2}\psi,\varphi\right)
 +i\left(\varphi^{\ast}\nabla\psi-\psi\nabla\varphi^{\ast},\mathbf{A}\right)}\\[2mm]
 {\displaystyle \quad \leq \|\nabla\psi\|_{\mathbf{L}^2}\|\nabla\varphi\|_{\mathbf{L}^2}+\|\mathbf{A}\|^{2}_{\mathbf{L}^6}\|\psi\|_{\mathcal{L}^6}\|\varphi\|_{\mathcal{L}^2}
 +\|\mathbf{A}\|_{\mathbf{L}^6}\big(\|\psi\|_{\mathcal{L}^3}\|\nabla\varphi\|_{\mathbf{L}^2}}\\[2mm]
 {\displaystyle \quad+\|\nabla\psi\|_{\mathbf{L}^2}\|\varphi\|_{\mathcal{L}^3}\big)\leq C\|\nabla\psi\|_{\mathbf{L}^2}
 \|\nabla\varphi\|_{\mathbf{L}^2},\quad \forall \mathbf{A}\in \mathbf{L}^6(\Omega),
 \,\,\,\psi,\varphi\in\mathcal{H}_{0}^{1}(\Omega),}
 \end{array}
 \end{equation}
 and
 \begin{equation}\label{eq:3-8}
 \begin{array}{@{}l@{}}
 {\displaystyle Q_3^{k}( \overline{\theta}_{\psi}^{k})= B(\mathbf{A}^{k-\frac{1}{2}};(\overline{\psi}^{k} - R_h\overline{\psi}^{k}),\overline{\theta}_{\psi}^{k})
 +B(\mathbf{A}^{k-\frac{1}{2}};( \psi^{k-\frac{1}{2}} - \overline{\psi}^{k}),\overline{\theta}_{\psi}^{k}).}
 \end{array}
 \end{equation}
Hence we obtain
 \begin{equation}\label{eq:3-9}
 {\displaystyle |Q_3^{k}(\overline{\theta}_{\psi}^{k})|\leq C\big\{h^{2r}+(\Delta t)^{4}\big\}+C\|\nabla\overline{\theta}_{\psi}^{k}\|^{2}_{\mathbf{L}^2}.}
 \end{equation}

To estimate $Q_4^{k}( \overline{\theta}_{\psi}^{k})$, we rewrite it as follows:
\begin{equation}\label{eq:3-9-0}
\begin{array}{@{}l@{}}
{\displaystyle Q_4^{k}( \overline{\theta}_{\psi}^{k}) = \left((\psi^{k-\frac12}-R_h\psi^{k-\frac12})\phi^{k-\frac12},\,\overline{\theta}_{\psi}^{k}\right)
 +\left( R_h(\psi^{k-\frac12}-\overline{\psi}^{k})\phi^{k-\frac12}, \,\overline{\theta}_{\psi}^{k}\right) }\\[2mm]
{\displaystyle \qquad \quad + \left( (R_h\overline{\psi}^{k}-\overline{\psi}^{k}_{h})\phi^{k-\frac12}, \,\overline{\theta}_{\psi}^{k}\right)+\left(\overline{\psi}^{k}_{h}(\phi^{k-\frac12}-I_{h}\phi^{k-\frac12}), \,\overline{\theta}_{\psi}^{k}\right)}\\[2mm]
{\displaystyle \qquad \quad +\left(\overline{\psi}^{k}_{h}I_{h}(\phi^{k-\frac12}-\overline{\phi}^{k}), \,\overline{\theta}_{\psi}^{k}\right)
 + \left( \overline{\psi}^{k}_{h} (I_{h}\overline{\phi}^{k}-\overline{\phi}^{k}_{h}), \,\overline{\theta}_{\psi}^{k}\right). }
\end{array}
\end{equation}

It follows from Lemma~\ref{lem3-1}, the regularity assumption (\ref{eq:2-10}), the properties of the interpolation operators and (\ref{eq:3-9-0}) that
\begin{equation}\label{eq:3-9-1}
\begin{array}{@{}l@{}}
{\displaystyle \vert Q_4^{k}( \overline{\theta}_{\psi}^{k}) \vert \leq C\big\{h^{2r}+(\Delta t)^{4}\big\}+C\left\{\|\nabla\overline{\theta}_{\psi}^{k}\|^{2}_{\mathbf{L}^2}
 +\|\nabla\overline{\theta}_{\phi}^{k}\|^{2}_{\mathbf{L}^2}\right\}.}
\end{array}
\end{equation}

We observe that
 \begin{equation}\label{eq:3-10}
 \begin{array}{@{}l@{}}
 {\displaystyle Q_5^{k}( \overline{\theta}_{\psi}^{k}) = \Big[B(\overline{\mathbf{A}}^{k}_{h}; R_h\overline{\psi}^{k},\overline{\theta}_{\psi}^{k})
 -B({\pi}_{h}\overline{\mathbf{A}}^{k}; R_h\overline{\psi}^{k},\overline{\theta}_{\psi}^{k})\Big]}\\[2mm]
 {\displaystyle \quad
 +\Big[B({\pi}_{h}\overline{\mathbf{A}}^{k}; R_h\overline{\psi}^{k},\overline{\theta}_{\psi}^{k})-B(\overline{\mathbf{A}}^{k}; R_h\overline{\psi}^{k},\overline{\theta}_{\psi}^{k})\Big]}\\[2mm]
 {\displaystyle \quad+\Big[B(\overline{\mathbf{A}}^{k}; R_h\overline{\psi}^{k},\overline{\theta}_{\psi}^{k})
 -B(\mathbf{A}^{k-\frac{1}{2}}; R_h\overline{\psi}^{k},\overline{\theta}_{\psi}^{k})\Big]}\\[2mm]
 {\displaystyle \quad \stackrel{\mathrm{def}}{=}Q_5^{k,1}( \overline{\theta}_{\psi}^{k})+ Q_5^{k,2}( \overline{\theta}_{\psi}^{k})+Q_5^{k,3}( \overline{\theta}_{\psi}^{k}).}
 \end{array}
 \end{equation}

It follows from Lemma~\ref{lem3-2} that
 \begin{equation}\label{eq:3-11}
 \begin{array}{@{}l@{}}
  {\displaystyle Q_5^{k,1}( \overline{\theta}_{\psi}^{k})=\left(R_h\overline{\psi}^{k}(\overline{\theta}_{\psi}^{k})^{\ast}(\overline{\mathbf{A}}^{k}_{h}
  +{\pi}_h\overline{\mathbf{A}}^{k}),\;\overline{\theta}^{k}_{\mathbf{A}}\right)}\\[2mm]
 {\displaystyle\qquad\quad+\mathrm{i}\left((\overline{\theta}_{\psi}^{k})^{\ast}\nabla R_h\overline{\psi}^{k}-R_h\overline{\psi}^{k}\nabla(\overline{\theta}_{\psi}^{k})^{\ast},
 \;\overline{\theta}^{k}_{\mathbf{A}}\right),}\\[2mm]
   {\displaystyle Q_5^{k,2}( \overline{\theta}_{\psi}^{k})=\left(R_h\overline{\psi}^{k}(\overline{\theta}_{\psi}^{k})^{\ast}({\pi}_h\overline{\mathbf{A}}^{k}+\overline{\mathbf{A}}^{k}),\; {\pi}_h\overline{\mathbf{A}}^{k}
   -\overline{\mathbf{A}}^{k}\right)}\\[2mm]
   {\displaystyle \qquad\quad+\mathrm{i}\left((\overline{\theta}_{\psi}^{k})^{\ast}\nabla R_h\overline{\psi}^{k}-R_h\overline{\psi}^{k}\nabla(\overline{\theta}_{\psi}^{k})^{\ast},\; {\pi}_h\overline{\mathbf{A}}^{k}-\overline{\mathbf{A}}^{k}\right),}\\[2mm]
  {\displaystyle Q_5^{k,3}( \overline{\theta}_{\psi}^{k})=\left(R_h\overline{\psi}^{k}(\overline{\theta}_{\psi}^{k})^{\ast}(\overline{\mathbf{A}}^{k}
  +\mathbf{A}^{k-\frac{1}{2}}),\;\overline{\mathbf{A}}^{k}-\mathbf{A}^{k-\frac{1}{2}}\right)}\\[2mm]
 {\displaystyle \qquad\quad +\mathrm{i}\left((\overline{\theta}_{\psi}^{k})^{\ast}\nabla R_h\overline{\psi}^{k}-R_h\overline{\psi}^{k}\nabla(\overline{\theta}_{\psi}^{k})^{\ast},\;
 \overline{\mathbf{A}}^{k}-\mathbf{A}^{k-\frac{1}{2}}\right).}
 \end{array}
 \end{equation}

Using (\ref{eq:3-84-0}), we prove
 \begin{equation}\label{eq:3-12}
 \begin{array}{@{}l@{}}
 {\displaystyle |Q_5^{k,1}( \overline{\theta}_{\psi}^{k})|\leq C\|\overline{\theta}^{k}_{\mathbf{A}}\|_{\mathbf{L}^6}\left\{\|\overline{\theta}_{\psi}^{k}\|_{\mathcal{L}^2}
 +\|\overline{\theta}_{\psi}^{k}\|_{\mathcal{L}^3}
+\|\nabla\overline{\theta}_{\psi}^{k}\|_{\mathbf{L}^2}\right\}}\\[2mm]
{\displaystyle\qquad \qquad\quad
    \leq C\left\{D(\overline{\theta}^{k}_{\mathbf{A}},\overline{\theta}^{k}_{\mathbf{A}}) + \|\nabla\overline{\theta}_{\psi}^{k}\|_{\mathbf{L}^2}^{2} \right\},  }\\[2mm]
 {\displaystyle |Q_5^{k,2}( \overline{\theta}_{\psi}^{k})|\leq C h^{r}\left\{\|\overline{\theta}_{\psi}^{k}\|_{\mathcal{L}^2}
 +\|\overline{\theta}_{\psi}^{k}\|_{\mathcal{L}^3}
+\|\nabla\overline{\theta}_{\psi}^{k}\|_{\mathbf{L}^2}\right\} \leq C \left\{h^{2r} +\|\nabla\overline{\theta}_{\psi}^{k}\|_{\mathbf{L}^2}^{2} \right\},}\\[2mm]
 {\displaystyle |Q_5^{k,3}( \overline{\theta}_{\psi}^{k})|\leq C(\Delta t)^{2}\left\{\|\overline{\theta}_{\psi}^{k}\|_{\mathcal{L}^2}
 +\|\overline{\theta}_{\psi}^{k}\|_{\mathcal{L}^3}
+\|\nabla\overline{\theta}_{\psi}^{k}\|_{\mathbf{L}^2}\right\} }\\[2mm]
{\displaystyle\qquad \qquad\quad\leq C\left\{(\Delta t)^{4} + \|\nabla\overline{\theta}_{\psi}^{k}\|_{\mathbf{L}^2}^{2} \right\},}
\end{array}
\end{equation}
and thus
\begin{equation}\label{eq:3-13}
 \begin{array}{@{}l@{}}
{\displaystyle |Q_5^{k}( \overline{\theta}_{\psi}^{k})|\leq C\big\{h^{2r}+(\Delta t)^{4}\big\} + C\left\{ D(\overline{\theta}^{k}_{\mathbf{A}},\overline{\theta}^{k}_{\mathbf{A}})+\|\nabla\overline{\theta}_{\psi}^{k}\|_{\mathbf{L}^2}^{2} \right\}. }
\end{array}
\end{equation}

From (\ref{eq:3-5-0}), summing over $k = 1,2,\cdots,m $ and combining (\ref{eq:3-6}), (\ref{eq:3-9}), (\ref{eq:3-9-1}) and
(\ref{eq:3-13}), we have
\begin{equation}\label{eq:3-15}
\begin{array}{@{}l@{}}
 {\displaystyle \|\theta^{m}_{\psi}\|^{2}_{\mathcal{L}^2}\leq C\big\{h^{2r}+(\Delta t)^{4}\big\}+C \Delta t
 \sum_{k=1}^{m}\left\{D(\overline{\theta}^{k}_{\mathbf{A}},\overline{\theta}^{k}_{\mathbf{A}}) + \|\nabla\overline{\theta}_{\psi}^{k}\|_{\mathbf{L}^2}^{2} + \|\nabla\overline{\theta}_{\phi}^{k}\|^{2}_{\mathbf{L}^2}\right\}}\\[2mm]
 {\displaystyle \quad \leq C\big\{h^{2r}+(\Delta t)^{4}\big\}
 +C\Delta t\sum_{k=0}^{m}\left\{D({\theta}^{k}_{\mathbf{A}},
 {\theta}^{k}_{\mathbf{A}}) + \|\nabla\theta_{\psi}^{k}\|^{2}_{\mathbf{L}^2} + \|\nabla {\theta}_{\phi}^{k}\|^{2}_{\mathbf{L}^2}\right\}.}
 \end{array}
 \end{equation}

To proceed further, we take $\varphi=\partial{\theta_{\psi}^{k}}=\frac{\displaystyle 1}{\displaystyle \Delta t}(\theta_{\psi}^{k}-\theta_{\psi}^{k-1})$, $ k=1,2,\cdots, m $ in (\ref{eq:3-3}),  to find
\begin{equation}\label{eq:3-16}
 {\displaystyle -2\mathrm{i}(\partial \theta^{k}_{\psi},\partial{\theta_{\psi}^{k}})+ B\left(\overline{\mathbf{A}}^{k}_{h};\overline{\theta}_{\psi}^{k},\partial{\theta_{\psi}^{k}}\right)
 =\sum_{j=1}^{5}Q_{j}^{k}(\partial{\theta_{\psi}^{k}}) .}
 \end{equation}

We take the real part of (\ref{eq:3-16}) and use (\ref{eq:3-2-5}) to get
\begin{equation}\label{eq:3-16-0}
\begin{array}{@{}l@{}}
{\displaystyle \frac{1}{2\Delta t}\left(B(\overline{\mathbf{A}}_{h}^{k};\theta_{\psi}^{k},\theta_{\psi}^{k})- B(\overline{\mathbf{A}}_{h}^{k-1};\theta_{\psi}^{k-1},\theta_{\psi}^{k-1}) \right) = \sum_{j=1}^{5}\mathrm{Re}\big[Q_{j}^{k}(\partial{\theta_{\psi}^{k}})\big]  }\\[2mm]
{\displaystyle \qquad + \left(\frac{1}{2}(\overline{\mathbf{A}}_{h}^{k}+\overline{\mathbf{A}}_{h}^{k-1})|\theta_{\psi}^{k-1}|^{2},\frac{1}{2}(\partial \mathbf{A}_{h}^{k}+\partial \mathbf{A}_{h}^{k-1})\right)        }\\[2mm]
{\displaystyle\qquad +\left(f(\theta_{\psi}^{k-1},\theta_{\psi}^{k-1}),\frac{1}{2}(\partial \mathbf{A}_{h}^{k}+\partial \mathbf{A}_{h}^{k-1})\right) ,}
\end{array}
\end{equation}
which leads to
\begin{equation}\label{eq:3-16-1}
\begin{array}{@{}l@{}}
{\displaystyle \frac{1}{2}B(\overline{\mathbf{A}}_{h}^{m};\theta_{\psi}^{m},\theta_{\psi}^{m}) = \frac{1}{2}B(\overline{\mathbf{A}}_{h}^{0};\theta_{\psi}^{0},\theta_{\psi}^{0}) + \Delta t \sum_{j=1}^{5} \sum_{k=1}^{m}\mathrm{Re}\big[Q_{j}^{k}(\partial{\theta_{\psi}^{k}})\big] }\\[2mm]
{\displaystyle \quad + \Delta t\sum_{k=1}^{m}\left(\frac{1}{2}(\overline{\mathbf{A}}_{h}^{k}+\overline{\mathbf{A}}_{h}^{k-1})|\theta_{\psi}^{k-1}|^{2},\,\overline{\partial \mathbf{A}}_{h}^{k}\right)}\\[2mm]
 {\displaystyle \quad +\Delta t\sum_{k=1}^{m} \left(f(\theta_{\psi}^{k-1},\theta_{\psi}^{k-1}),\,\overline{\partial \mathbf{A}}_{h}^{k}\right).}
\end{array}
\end{equation}

It follows from (\ref{eq:3-84-0}) that
\begin{equation}\label{eq:3-16-2}
\begin{array}{@{}l@{}}
{\displaystyle   \sum_{k=1}^{m}\left(\frac{1}{2}(\overline{\mathbf{A}}_{h}^{k}+\overline{\mathbf{A}}_{h}^{k-1})|\theta_{\psi}^{k-1}|^{2},\,\overline{\partial \mathbf{A}}_{h}^{k}\right) \leq \sum_{k=1}^{m} \| \overline{\mathbf{A}}_{h}^{k}+\overline{\mathbf{A}}_{h}^{k-1}\|_{\mathbf{L}^{6}} \|\theta_{\psi}^{k-1}\|^{2}_{\mathcal{L}^{6}} \|\overline{\partial \mathbf{A}}_{h}^{k}\|_{\mathbf{L}^{2}}  }\\[2mm]
{\displaystyle \qquad  \leq C\sum_{k=1}^{m}\|\theta_{\psi}^{k-1}\|^{2}_{\mathcal{L}^{6}}   \leq C\sum_{k=0}^{m}\|\nabla\theta_{\psi}^{k}\|^{2}_{\mathbf{L}^{2}}. }\\[2mm]
\end{array}
\end{equation}

Combining (\ref{eq:3-81}) and (\ref{eq:3-82}) gives
\begin{equation}\label{eq:3-16-3}
\begin{array}{@{}l@{}}
{\displaystyle \sum_{k=1}^{m}\left(f(\theta_{\psi}^{k-1},\theta_{\psi}^{k-1}),\;\overline{\partial \mathbf{A}}_{h}^{k}\right) =  \sum_{k=1}^{m}\left(f(\theta_{\psi}^{k-1},\theta_{\psi}^{k-1}),\;\overline{\partial \theta}_{\mathbf{A}}^{k}\right) + \sum_{k=1}^{m}\left(f(\theta_{\psi}^{k-1},\theta_{\psi}^{k-1}),\;\overline{\partial {\pi}_{h}\mathbf{A}}^{k}\right)   }\\[2mm]
{\displaystyle \quad\leq \sum_{k=1}^{m}\left(f(\theta_{\psi}^{k-1},\theta_{\psi}^{k-1}),\;\overline{\partial \theta}_{\mathbf{A}}^{k}\right) +  C\sum_{k=1}^{m}\|\nabla \theta_{\psi}^{k-1}\|_{\mathbf{L}^{2}} \| \theta_{\psi}^{k-1}\|_{\mathcal{L}^{6}}  \|\overline{\partial {\pi}_{h}\mathbf{A}}^{k}\|_{\mathbf{L}^{3}}  }\\[2mm]
{\displaystyle \quad\leq Ch^{2r}+C \sum_{k=0}^{m}\left\{ \|\nabla \theta_{\psi}^{k}\|^{2}_{\mathbf{L}^{2}} + D(\theta_{\mathbf{A}}^{k},\theta_{\mathbf{A}}^{k})+ \|\partial \theta_{\mathbf{A}}^{k} \|^{2}_{\mathbf{L}^{2}}\right\}.}
\end{array}
\end{equation}

Substituting (\ref{eq:3-16-2}) and (\ref{eq:3-16-3}) into (\ref{eq:3-16-1}), we have
\begin{equation}\label{eq:3-16-4}
\begin{array}{@{}l@{}}
{\displaystyle \frac{1}{2}B(\overline{\mathbf{A}}_{h}^{m};\theta_{\psi}^{m},\theta_{\psi}^{m}) \leq
 C h^{2r} + \Delta t \sum_{j=1}^{5} \sum_{k=1}^{m}\mathrm{Re}\big[Q_{j}^{k}(\partial{\theta_{\psi}^{k}})\big] }\\[2mm]
 {\displaystyle \qquad +C \Delta t\sum_{k=0}^{m}\left( \|\nabla \theta_{\psi}^{k}\|^{2}_{\mathbf{L}^{2}} + D(\theta_{\mathbf{A}}^{k},\theta_{\mathbf{A}}^{k}) + \|\partial \theta_{\mathbf{A}}^{k} \|^{2}_{\mathbf{L}^{2}}\right)  .}
 \end{array}
\end{equation}

We now proceed to estimate $\sum_{k=1}^{m}\mathrm{Re}\big[Q_{j}^{k}(\partial{\theta_{\psi}^{k}})\big] $, $j=1,\cdots,5$.
By virtue of (\ref{eq:3-2-2}), we get
\begin{equation}\label{eq:3-17}
\begin{array}{@{}l@{}}
{\displaystyle \Delta t\sum_{k=1}^{m} Q_{1}^{k}(\partial{\theta_{\psi}^{k}})=2\mathrm{i}\sum_{k=1}^{m}\left(\partial R_h\psi^{k}-\frac{\partial \psi^{k-\frac{1}{2}}}{\partial t},\theta_{\psi}^{k}-{\theta_{\psi}^{k-1}}\right)}\\[2mm]
{\displaystyle \quad= 2\mathrm{i}\left(\partial R_h\psi^{m}-\frac{\partial \psi^{m-\frac{1}{2}}}{\partial t},\theta_{\psi}^{m}\right)-2\mathrm{i}\left(\partial R_h\psi^{1}
-\frac{\partial \psi^{\frac{1}{2}}}{\partial t},\theta_{\psi}^{0}\right)}\\[2mm]
{\displaystyle \quad-2\mathrm{i}\sum_{k=1}^{m-1}\left(\partial R_h\psi^{k+1}
-\partial R_h\psi^{k}-\frac{\partial \psi^{k+\frac{1}{2}}}{\partial t}
+\frac{\partial \psi^{k-\frac{1}{2}}}{\partial t},\theta_{\psi}^{k}\right).}
\end{array}
\end{equation}

From (\ref{eq:3-2-3}) and (\ref{eq:3-17}), we deduce
\begin{equation}\label{eq:3-19}
\begin{array}{@{}l@{}}
{\displaystyle |\Delta t\sum_{k=1}^{m} Q_{1}^{k}(\partial{\theta_{\psi}^{k}})|\leq C\big\{h^{2r+2}+(\Delta t)^{4}\big\}
+ C\|\theta_{\psi}^{m}\|_{\mathcal{L}^2}^{2}+C\Delta t \sum_{k=1}^{m-1}{\|\theta_{\psi}^{k}\|_{\mathcal{L}^2}^{2}}.}
\end{array}
\end{equation}

We observe that
\begin{equation*}
\begin{array}{@{}l@{}}
{\displaystyle \Delta t Q_{2}^{k}(\partial{\theta_{\psi}^{k}}) = 2 V_0\left(\psi^{k-\frac{1}{2}}-\overline{\psi}_{h}^{k},\theta_{\psi}^{k}-\theta_{\psi}^{k-1}\right) }\\[2mm]
{\displaystyle  \quad= 2V_0\left(\psi^{k-\frac{1}{2}}
-R_h\overline{\psi}^{k},\theta_{\psi}^{k}-\theta_{\psi}^{k-1}\right)-2V_0\left(\frac{1}{2}(\theta_{\psi}^{k}
+\theta_{\psi}^{k-1}),\theta_{\psi}^{k}-\theta_{\psi}^{k-1}\right) }\\[2mm]
{\displaystyle \quad \stackrel{\mathrm{def}}{=}J_2^{k,1}+J_2^{k,2}.}
\end{array}
\end{equation*}

It is not difficult to prove
\begin{equation}\label{eq:3-20}
\begin{array}{@{}l@{}}
{\displaystyle |\sum_{k=1}^{m} J_2^{k,1}| \leq C\big\{h^{2r+2}+(\Delta t)^{4}\big\}
 +C\|\theta_{\psi}^{m}\|_{\mathcal{L}^2}^{2}+ C\Delta t \sum_{k=1}^{m-1}{\|\theta_{\psi}^{k}\|_{\mathcal{L}^2}^{2}}.}
\end{array}
\end{equation}

We estimate $|\sum\limits_{k=1}^{m} \mathrm{Re}[J_2^{k,2}]|$ by
\begin{equation}\label{eq:3-21}
\begin{array}{@{}l@{}}
{\displaystyle |\sum\limits_{k=1}^{m} \mathrm{Re}[J_2^{k,2}]|=|-V_0(\|\theta_{\psi}^{m}\|_{\mathcal{L}^2}^{2}-\|\theta_{\psi}^{0}\|_{\mathcal{L}^2}^{2})|
\leq C\|\theta_{\psi}^{m}\|_{\mathcal{L}^2}^{2}+C h^{2r+2}.}
\end{array}
\end{equation}

Hence we have
\begin{equation}\label{eq:3-22}
\begin{array}{@{}l@{}}
{\displaystyle  |\Delta t\sum_{k=1}^{m} \mathrm{Re}\big[ Q_{2}^{k}(\partial{\theta_{\psi}^{k}})\big]|\leq  |\sum\limits_{k=1}^{m} \mathrm{Re}[J_2^{k,1}]|
+ |\sum\limits_{k=1}^{m} \mathrm{Re}[J_2^{k,2}]|\leq |\sum_{k=1}^{m}J_2^{k,1}|}\\[2mm]
{\displaystyle \quad + |\sum\limits_{k=1}^{m} \mathrm{Re}[J_2^{k,2}]|\leq C\big\{h^{2r+2}+(\Delta t)^{4}\big\}
+C \|\theta_{\psi}^{m}\|_{\mathcal{L}^2}^{2}+ C\Delta t \sum_{k=1}^{m-1}{\|\theta_{\psi}^{k}\|_{\mathcal{L}^2}^{2}}.}
\end{array}
\end{equation}

From the definition of the bilinear functional $B(\mathbf{A};\psi,\varphi)$ in (\ref{eq:2-7}), we rewrite $\Delta t  Q_{3}^{k}(\partial{\theta_{\psi}^{k}})$ as follows:
\begin{equation}\label{eq:3-23}
\begin{array}{@{}l@{}}
{\displaystyle \Delta t  Q_{3}^{k}(\partial{\theta_{\psi}^{k}})= \left(\nabla( \psi^{k-\frac{1}{2}}-R_h\overline{\psi}^{k}),\;\nabla
(\theta_{\psi}^{k}-\theta_{\psi}^{k-1})\right)}\\[2mm]
{\displaystyle\quad \qquad +\left(|\mathbf{A}^{k-\frac{1}{2}}|^2(\psi^{k-\frac{1}{2}}-R_h\overline{\psi}^{k}),\theta_{\psi}^{k}-\theta_{\psi}^{k-1}\right)}\\[2mm]
{\displaystyle \quad \qquad+ \mathrm{i}\left(\nabla(\psi^{k-\frac{1}{2}}-R_h\overline{\psi}^{k})\mathbf{A}^{k-\frac{1}{2}},\;\theta_{\psi}^{k}
-\theta_{\psi}^{k-1}\right)}\\[2mm]
{\displaystyle \quad \qquad-\mathrm{i}\left((\psi^{k-\frac{1}{2}}-R_h\overline{\psi}^{k})\mathbf{A}^{k-\frac{1}{2}},\;\nabla \theta_{\psi}^{k}
-\nabla \theta_{\psi}^{k-1}\right).}
\end{array}
\end{equation}

Analogous to the estimate of $\sum_{k=1}^{m}  Q_{1}^{k}(\partial{\theta_{\psi}^{k}})$, we use (\ref{eq:3-2-2}), (\ref{eq:3-2-3}), the regularity assumption (\ref{eq:2-10}) and Young's inequality to prove
\begin{equation}\label{eq:3-23-0}
{\displaystyle  |\Delta t\sum_{k=1}^{m}  Q_{3}^{k}(\partial{\theta_{\psi}^{k}})| \leq C\{h^{2r}+(\Delta t)^{4}\} +C \Vert\theta_{\psi}^{m}\Vert_{\mathcal{L}^2}^{2}
+ \frac{1}{16} \Vert\nabla\theta_{\psi}^{m}\Vert_{\mathbf{L}^2}^{2}+C\Delta t \sum_{k=0}^{m}{\Vert\nabla\theta_{\psi}^{k}\Vert_{\mathbf{L}^2}^{2}}.}
\end{equation}
Due to space limitations, we omit the proof of (\ref{eq:3-23-0}).

We observe that
\begin{equation*}
{\displaystyle \phi^{k-\frac12}\psi^{k-\frac12}-\overline{\phi}^{k}_{h}\overline{\psi}^{k}_{h} = \phi^{k-\frac12}\psi^{k-\frac12}-I_{h}\overline{\phi}^{k}R_h\overline{\psi}^{k}
-R_h\overline{\psi}^{k}\overline{\theta}_{\phi}^{k}-\overline{\phi}^{k}_{h}\overline{\theta}_{\psi}^{k},}
\end{equation*}
and have
\begin{equation}\label{eq:3-24}
 \begin{array}{@{}l@{}}
 {\displaystyle \Delta t Q_{4}^{k}(\partial{\theta_{\psi}^{k}}) = \left(\phi^{k-\frac12}\psi^{k-\frac12}-I_{h}\overline{\phi}^{k}R_h\overline{\psi}^{k}, \;\theta_{\psi}^{k}
 -\theta_{\psi}^{k-1} \right)}\\[2mm]
 {\displaystyle \qquad - \left( R_h\overline{\psi}^{k}\overline{\theta}_{\phi}^{k}, \;\theta_{\psi}^{k}-\theta_{\psi}^{k-1} \right) - \left(\overline{\phi}^{k}_{h}\overline{\theta}_{\psi}^{k}, \;\theta_{\psi}^{k}-\theta_{\psi}^{k-1} \right) }.
\end{array}
\end{equation}

The first two terms can be estimated as follows.
\begin{equation}\label{eq:4-25}
 \begin{array}{@{}l@{}}
 {\displaystyle |\sum_{k=1}^{m}\left(\phi^{k-\frac12}\psi^{k-\frac12}-I_{h}\overline{\phi}^{k}R_h\overline{\psi}^{k}, \;\theta_{\psi}^{k}-\theta_{\psi}^{k-1} \right)| \leq C\{h^{2r}+(\Delta t)^{4}\}  }\\[2mm]
 {\displaystyle \qquad\quad +C \Vert\theta_{\psi}^{m}\Vert_{\mathcal{L}^2}^{2} + C\Delta t \sum_{k=0}^{m}{\Vert\theta_{\psi}^{k}\Vert_{\mathcal{L}^2}^{2}}, }\\[2mm]
 {\displaystyle  |\sum_{k=1}^{m} \left( R_h\overline{\psi}^{k}\overline{\theta}_{\phi}^{k}, \;\theta_{\psi}^{k}-\theta_{\psi}^{k-1} \right) |  \leq Ch^{2r} + C \Vert\theta_{\psi}^{m}\Vert_{\mathcal{L}^2}^{2} +\frac{1}{8}\Vert\nabla\overline{\theta}_{\phi}^{m}\Vert_{\mathbf{L}^2}^{2}  }\\[2mm]
 {\displaystyle \qquad\quad+ C\Delta t \sum_{k=0}^{m}\left(\Vert\theta_{\psi}^{k}\Vert_{\mathcal{L}^2}^{2}+   \Vert\nabla\theta_{\phi}^{k}\Vert_{\mathbf{L}^2}^{2}+\Vert\partial \theta_{\phi}^{k}\Vert_{{L}^2}^{2} \right). }
\end{array}
\end{equation}

By applying (\ref{eq:3-93}), the last term on the right hand side of (\ref{eq:3-24}) can be bounded as follows:
\begin{equation}\label{eq:3-26}
 \begin{array}{@{}l@{}}
 {\displaystyle  \sum_{k=1}^{m} \mathrm{Re}\big[\left(\overline{\phi}^{k}_{h}\overline{\theta}_{\psi}^{k}, \;\theta_{\psi}^{k}-\theta_{\psi}^{k-1} \right)  \big]=\frac{1}{2} \sum_{k=1}^{m} \left(\overline{\phi}^{k}_{h}, \;\vert \theta_{\psi}^{k} \vert^{2} -\vert \theta_{\psi}^{k-1}\vert^{2} \right) }\\[2mm]
 {\displaystyle = \frac{1}{2}\left( \overline{\phi}^{m}_{h},\;\vert \theta_{\psi}^{m} \vert^{2}\right) -   \frac{1}{2}\left( \overline{\phi}^{0}_{h},\;\vert \theta_{\psi}^{0} \vert^{2}\right) - \frac{1}{2} \sum_{k=1}^{m}\left( \overline{\phi}^{k}_{h} - \overline{\phi}^{k-1}_{h},  \vert \theta_{\psi}^{k-1} \vert^{2}  \right)} \\[2mm]
 {\displaystyle \leq \frac{1}{2} \|\overline{\phi}^{m}_{h}\|_{L^6} \|\theta_{\psi}^{m}\|_{\mathcal{L}^{3}} \|\theta_{\psi}^{m}\|_{\mathcal{L}^{2}}  + Ch^{2r} + C\Delta t \sum_{k=1}^{m} \|\partial \overline{\phi}^{k}_{h} \|_{L^2}\| \theta_{\psi}^{k-1}\|_{\mathcal{L}^{4}}^2 }\\[2mm]
 {\displaystyle \leq C\|\theta_{\psi}^{m}\|_{\mathcal{L}^{3}} \|\theta_{\psi}^{m}\|_{\mathcal{L}^{2}}
 + Ch^{2r} + C\Delta t \sum_{k=1}^{m}\|\theta_{\psi}^{k-1}\|_{\mathcal{L}^{4}}^2 }\\[2mm]
 {\displaystyle  \leq Ch^{2r} + C\|\theta_{\psi}^{m}\|^{2}_{\mathcal{L}^{2}} + \frac{1}{16} \|\nabla\theta_{\psi}^{m}\|^{2}_{\mathbf{L}^{2}} +C\Delta t \sum_{k=0}^{m} \|\nabla \theta_{\psi}^{k}\|_{\mathbf{L}^{2}}^2.}
 \end{array}
\end{equation}

Combining (\ref{eq:3-24})-(\ref{eq:3-26}) gives
\begin{equation}\label{eq:3-27}
  \begin{array}{@{}l@{}}
 {\displaystyle  \Delta t \sum_{k=1}^{m} \mathrm{Re}\big[Q_{4}^{k}(\partial{\theta_{\psi}^{k}}) \big] \leq  C\{h^{2r}+(\Delta t)^{4}\} + C \Vert\theta_{\psi}^{m}\Vert_{\mathcal{L}^2}^{2} +\frac{1}{8}\Vert\nabla\overline{\theta}_{\phi}^{m}\Vert_{\mathbf{L}^2}^{2} + \frac{1}{16} \|\nabla\theta_{\psi}^{m}\|^{2}_{\mathbf{L}^{2}}  }\\[2mm]
  {\displaystyle \qquad\quad+ C\Delta t \sum_{k=0}^{m}\left(\Vert\nabla \theta_{\psi}^{k}\Vert_{\mathbf{L}^2}^{2}+   \Vert\nabla\theta_{\phi}^{k}\Vert_{\mathbf{L}^2}^{2} + \Vert\partial \theta_{\phi}^{k}\Vert_{{L}^2}^{2} \right). }
\end{array}
\end{equation}

 $\Delta t Q_{5}^{k}(\partial{\theta_{\psi}^{k}})$ can be decomposed as follows:
\begin{equation}\label{eq:3-30}
\begin{array}{@{}l@{}}
{\displaystyle   \Delta t Q_{5}^{k}(\partial{\theta_{\psi}^{k}})= \Big[B(\mathbf{A}^{k-\frac{1}{2}}; R_h\overline{\psi}^{k},\theta_{\psi}^{k}-\theta_{\psi}^{k-1})
-B(\overline{\mathbf{A}}^{k}; R_h\overline{\psi}^{k},\theta_{\psi}^{k}-\theta_{\psi}^{k-1})\Big]}\\[2mm]
{\displaystyle  \qquad\quad \qquad +\Big[B(\overline{\mathbf{A}}^{k}; R_h\overline{\psi}^{k},\theta_{\psi}^{k}-\theta_{\psi}^{k-1})-B({\pi}_{h}\overline{\mathbf{A}}^{k}; R_h\overline{\psi}^{k},\theta_{\psi}^{k}-\theta_{\psi}^{k-1})\Big]}\\[2mm]
{\displaystyle  \qquad \quad \qquad +\Big[B({\pi}_{h}\overline{\mathbf{A}}^{k}; R_h\overline{\psi}^{k},\theta_{\psi}^{k}-\theta_{\psi}^{k-1})-B(\overline{\mathbf{A}}^{k}_{h}; R_h\overline{\psi}^{k},\theta_{\psi}^{k}-\theta_{\psi}^{k-1})\Big]}\\[2mm]
{\displaystyle\qquad \qquad\quad \stackrel{\mathrm{def}}{=}R_5^{k,1}+R_5^{k,2}+R_5^{k,3}.}
\end{array}
\end{equation}

Following the lines of the proof of (\ref{eq:3-27}) and using (\ref{eq:3-2-5}), we prove
\begin{equation}\label{eq:3-34}
\begin{array}{@{}l@{}}
{\displaystyle  |\sum_{k=1}^{m} R_5^{k,1}| + |\sum_{k=1}^{m} R_5^{k,2}| \leq C\big\{h^{2r}+(\Delta t)^{4}\big\}
 +C \|\theta_{\psi}^{m}\|_{\mathcal{L}^2}^{2}}\\[2mm]
{\displaystyle \quad+ \frac{1}{16} \|\nabla\theta_{\psi}^{m}\|_{\mathbf{L}^2}^{2}+C\Delta t \sum_{k=1}^{m}\|\nabla\theta_{\psi}^{k}\|_{\mathbf{L}^2}^{2}.}
\end{array}
\end{equation}

To estimate $|\sum\limits_{k=1}^{m} R_5^{k,3}|$, we rewrite it as follows:
\begin{equation}\label{eq:3-35}
\begin{array}{@{}l@{}}
{\displaystyle  \sum_{k=1}^{m}R_5^{k,3}=\sum_{k=1}^{m}{ \left(R_h\overline{\psi}^{k}({\pi}_{h}\overline{\mathbf{A}}^{k}+\overline{\mathbf{A}}^{k}_{h})({\pi}_{h}\overline{\mathbf{A}}^{k}
-\overline{\mathbf{A}}^{k}_{h}),\;\theta_{\psi}^{k}-\theta_{\psi}^{k-1}\right)}}\\[2mm]
{\displaystyle\quad\quad\quad\quad-\sum_{k=1}^{m}{\mathrm{i}\left( R_h\overline{\psi}^{k}({\pi}_{h}\overline{\mathbf{A}}^{k}-\overline{\mathbf{A}}^{k}_{h}),\;
\nabla\theta_{\psi}^{k}-\nabla\theta_{\psi}^{k-1}\right)}}\\[2mm]
{\displaystyle \quad\quad\quad\quad+\sum_{k=1}^{m}{\mathrm{i}\left(\nabla R_h\overline{\psi}^{k}({\pi}_{h}\overline{\mathbf{A}}^{k}
-\overline{\mathbf{A}}^{k}_{h}),\;\theta_{\psi}^{k}-\theta_{\psi}^{k-1}\right)}}\\[2mm]
{\displaystyle\quad\quad\quad \stackrel{\mathrm{def}}{=}K_1+K_2+K_3.}
\end{array}
\end{equation}

Note that
\begin{equation}\label{eq:3-36}
\begin{array}{@{}l@{}}
{\displaystyle K_1= \sum_{k=1}^{m}\left(R_h\overline{\psi}^{k}({\pi}_{h}\overline{\mathbf{A}}^{k}+\overline{\mathbf{A}}^{k}_{h})
({\pi}_{h}\overline{\mathbf{A}}^{k}-\overline{\mathbf{A}}^{k}_{h}),\;\theta_{\psi}^{k}-\theta_{\psi}^{k-1}\right)}\\[2mm]
{\displaystyle \quad = -\left(R_h\overline{\psi}^{m}({\pi}_{h}\overline{\mathbf{A}}^{m}
+\overline{\mathbf{A}}^{m}_{h})
\overline{\theta}_{\mathbf{A}}^{m},\;\theta_{\psi}^{m}\right) + \left(R_h\overline{\psi}^{0}({\pi}_{h}\overline{\mathbf{A}}^{0}+\overline{\mathbf{A}}^{0}_{h})
\overline{\theta}_{\mathbf{A}}^{0},\;\theta_{\psi}^{0}\right)}\\[2mm]
{\displaystyle \quad + \sum_{k=1}^{m}\left(R_h\overline{\psi}^{k}({\pi}_{h}\overline{\mathbf{A}}^{k}
+\overline{\mathbf{A}}^{k}_{h})
\overline{\theta}_{\mathbf{A}}^{k}-R_h\overline{\psi}^{k-1}({\pi}_{h}\overline{\mathbf{A}}^{k-1}
+\overline{\mathbf{A}}^{k-1}_{h})
\overline{\theta}_{\mathbf{A}}^{k-1},\;\theta_{\psi}^{k-1}\right).}
\end{array}
\end{equation}

By applying Young's inequality and (\ref{eq:3-84-0}), we can estimate the first two terms on the right side of (\ref{eq:3-36}) by
\begin{equation}\label{eq:3-37}
\begin{array}{@{}l@{}}
{\displaystyle |\left(R_h\overline{\psi}^{m}({\pi}_{h}\overline{\mathbf{A}}^{m}
+\overline{\mathbf{A}}^{m}_{h})\overline{\theta}_{\mathbf{A}}^{m},\;\theta_{\psi}^{m}\right)|
+|\left(R_h\overline{\psi}^{0}({\pi}_{h}\overline{\mathbf{A}}^{0}
+\overline{\mathbf{A}}^{0}_{h})\overline{\theta}_{\mathbf{A}}^{0},\;\theta_{\psi}^{0}\right)|}\\[2mm]
{\displaystyle \quad\leq \|R_h\overline{\psi}^{m}\|_{\mathcal{L}^6}\| {\pi}_{h}\overline{\mathbf{A}}^{m}+\overline{\mathbf{A}}^{m}_{h}\|_{\mathbf{L}^6}
\|\overline{\theta}_{\mathbf{A}}^{m}\|_{\mathbf{L}^6}\|\theta_{\psi}^{m}\|_{\mathcal{L}^2} + Ch^{2r}}\\[2mm]
{\displaystyle \quad\leq C  \|\overline{\theta}_{\mathbf{A}}^{m}\|_{\mathbf{H}^1}\|\theta_{\psi}^{m}\|_{\mathcal{L}^2}
+ C h^{2r} \leq C D^{\frac{1}{2}}(\overline{\theta}_{\mathbf{A}}^{m},\overline{\theta}_{\mathbf{A}}^{m})
\|\theta_{\psi}^{m}\|_{\mathcal{L}^2}+ Ch^{2r}}\\[2mm]
{\displaystyle \quad\leq \frac{1}{32}D(\overline{\theta}_{\mathbf{A}}^{m},\overline{\theta}_{\mathbf{A}}^{m})+ C\|\theta_{\psi}^{m}\|_{\mathcal{L}^2}^{2}+ Ch^{2r}.}
\end{array}
\end{equation}

Observing
\begin{equation}\label{eq:3-38}
\begin{array}{@{}l@{}}
{\displaystyle \left(R_h\overline{\psi}^{k}({\pi}_{h}\overline{\mathbf{A}}^{k}+\overline{\mathbf{A}}^{k}_{h})
\overline{\theta}_{\mathbf{A}}^{k}-R_h\overline{\psi}^{k-1}({\pi}_{h}\overline{\mathbf{A}}^{k-1}
+\overline{\mathbf{A}}^{k-1}_{h})\overline{\theta}_{\mathbf{A}}^{k-1},\;\theta_{\psi}^{k-1}\right)}\\[2mm]
{\displaystyle \quad= \Delta t\left(R_h\overline{\psi}^{k}({\pi}_{h}\overline{\mathbf{A}}^{k}
+\overline{\mathbf{A}}^{k}_{h})\frac{\overline{\theta}_{\mathbf{A}}^{k}-\overline{\theta}_{\mathbf{A}}^{k-1}}{\Delta t},\;\theta_{\psi}^{k-1}\right)}\\[2mm]
{\displaystyle \quad+ \Delta t\left(\frac{R_h\overline{\psi}^{k}-R_h\overline{\psi}^{k-1}}{\Delta t}({\pi}_{h}\overline{\mathbf{A}}^{k}+\overline{\mathbf{A}}^{k}_{h})\overline{\theta}_{\mathbf{A}}^{k-1},\;\theta_{\psi}^{k-1}\right)}\\[2mm]
{\displaystyle \quad+ \Delta t \left(R_h\overline{\psi}^{k-1}\overline{\theta}_{\mathbf{A}}^{k-1}\left(\frac{{\pi}_{h}\overline{\mathbf{A}}^{k}
-{\pi}_{h}\overline{\mathbf{A}}^{k-1}}{\Delta t}+\frac{\overline{\mathbf{A}}_{h}^{k}-\overline{\mathbf{A}}_{h}^{k-1}}{\Delta t}\right),\;\theta_{\psi}^{k-1}\right),}
\end{array}
\end{equation}
and using (\ref{eq:3-84-0}), we get
\begin{equation}\label{eq:3-39}
\begin{array}{@{}l@{}}
{\displaystyle |\left(R_h\overline{\psi}^{k}({\pi}_{h}\overline{\mathbf{A}}^{k}
+\overline{\mathbf{A}}^{k}_{h})\overline{\theta}_{\mathbf{A}}^{k}-R_h\overline{\psi}^{k-1}
({\pi}_{h}\overline{\mathbf{A}}^{k-1}+\overline{\mathbf{A}}^{k-1}_{h})
\overline{\theta}_{\mathbf{A}}^{k-1},\;\theta_{\psi}^{k-1}\right)|}\\[2mm]
{\displaystyle \quad\leq \Delta t \|R_h\overline{\psi}^{k}\|_{\mathcal{L}^6}\| {\pi}_{h}\overline{\mathbf{A}}^{k}+\overline{\mathbf{A}}^{k}_{h}\|_{\mathbf{L}^6}
\|\frac{1}{\Delta t} (\overline{\theta}_{\mathbf{A}}^{k}-\overline{\theta}_{\mathbf{A}}^{k-1})\|_{\mathbf{L}^2}\|\theta_{\psi}^{k-1}\|_{\mathcal{L}^6}}\\[2mm]
{\displaystyle \quad + \Delta t \|\frac{R_h\overline{\psi}^{k}-R_h\overline{\psi}^{k-1}}{\Delta t}\|_{\mathcal{L}^2}\| {\pi}_{h}\overline{\mathbf{A}}^{k}+\overline{\mathbf{A}}^{k}_{h}\|_{\mathbf{L}^6}
\|\overline{\theta}_{\mathbf{A}}^{k-1}\|_{\mathbf{L}^6}\|\theta_{\psi}^{k-1}\|_{\mathcal{L}^6}}\\[2mm]
{\displaystyle \quad+\Delta t \|R_h\overline{\psi}^{k-1}\|_{\mathcal{L}^6}\| \overline{\theta}_{\mathbf{A}}^{k-1}\|_{\mathbf{L}^6}\|\frac{{\pi}_{h}\overline{\mathbf{A}}^{k}
-{\pi}_{h}\overline{\mathbf{A}}^{k-1}}{\Delta t}+\frac{\overline{\mathbf{A}}_{h}^{k}-\overline{\mathbf{A}}_{h}^{k-1}}{\Delta t}\|_{\mathbf{L}^2}\|\theta_{\psi}^{k-1}\|_{\mathcal{L}^6}}\\[2mm]
{\displaystyle \quad\leq C\Delta t\left(\|\partial \overline{\theta}_{\mathbf{A}}^{k}\|_{\mathbf{L}^2}
 + \|\overline{\theta}_{\mathbf{A}}^{k-1}\|_{\mathbf{H}^1}\right)\|\theta_{\psi}^{k-1}\|_{\mathcal{H}^1}}\\[2mm]
{\displaystyle \quad\leq C\Delta t\left(\|\partial \overline{\theta}_{\mathbf{A}}^{k}\|^{2}_{\mathbf{L}^2} + D(\overline{\theta}_{\mathbf{A}}^{k-1},\overline{\theta}_{\mathbf{A}}^{k-1}) + \|\nabla\theta_{\psi}^{k-1}\|_{\mathbf{L}^2}^{2} \right). }
\end{array}
\end{equation}

From (\ref{eq:3-36})-(\ref{eq:3-39}), we thus have
\begin{equation}\label{eq:3-40}
\begin{array}{@{}l@{}}
{\displaystyle |K_1|= |\sum_{k=1}^{m}\left(R_h\overline{\psi}^{k}({\pi}_{h}\overline{\mathbf{A}}^{k}
+\overline{\mathbf{A}}^{k}_{h})({\pi}_{h}\overline{\mathbf{A}}^{k}-\overline{\mathbf{A}}^{k}_{h}),\;\theta_{\psi}^{k}
-\theta_{\psi}^{k-1}\right)|}\\[2mm]
{\displaystyle \quad\leq \frac{1}{16}D(\overline{\theta}_{\mathbf{A}}^{m},\overline{\theta}_{\mathbf{A}}^{m})
+ C\|\theta_{\psi}^{m}\|_{\mathcal{L}^2}^{2}+ C h^{2r}}\\[2mm]
{\displaystyle\quad + C\Delta t \sum_{k=0}^{m}\left(\|\partial {\theta}_{\mathbf{A}}^{k}\|_{\mathbf{L}^2}^{2}
+ D({\theta}_{\mathbf{A}}^{k},{\theta}_{\mathbf{A}}^{k})+ \|\nabla\theta_{\psi}^{k}\|_{\mathbf{L}^2}^{2}\right)}\\[2mm]
\end{array}
\end{equation}

From (\ref{eq:3-35}), using (\ref{eq:3-2-2}) and integrating by parts, we get
\begin{equation}\label{eq:3-41}
\begin{array}{@{}l@{}}
{\displaystyle K_2=\sum_{k=1}^{m}\left( R_h\overline{\psi}^{k}({\pi}_{h}\overline{\mathbf{A}}^{k}
-\overline{\mathbf{A}}^{k}_{h}),\;\nabla\theta_{\psi}^{k}-\nabla\theta_{\psi}^{k-1}\right)}\\[2mm]
{\displaystyle\quad = -\left(R_h\overline{\psi}^{m}\overline{\theta}_{\mathbf{A}}^{m},\;\nabla\theta_{\psi}^{m}\right)
+\left(R_h\overline{\psi}^{0}\overline{\theta}_{\mathbf{A}}^{0},\;\nabla\theta_{\psi}^{0}\right)}\\[2mm]
{\displaystyle\quad+ \sum_{k=1}^{m}\left(R_h\overline{\psi}^{k}\overline{\theta}_{\mathbf{A}}^{k}- R_h\overline{\psi}^{k-1}\overline{\theta}_{\mathbf{A}}^{k-1},\;\nabla\theta_{\psi}^{k-1}\right)}\\[2mm]
{\displaystyle \quad=\left(\nabla R_h\overline{\psi}^{m}\overline{\theta}_{\mathbf{A}}^{m},\;\theta_{\psi}^{m}\right)
 + \left(R_h\overline{\psi}^{m}\nabla\cdot\overline{\theta}_{\mathbf{A}}^{m},\;\theta_{\psi}^{m}\right)}\\[2mm]
{\displaystyle \quad +\left(R_h\overline{\psi}^{0}\overline{\theta}_{\mathbf{A}}^{0},\;\nabla\theta_{\psi}^{0}\right)
+\sum_{k=1}^{m}\left(R_h\overline{\psi}^{k}\overline{\theta}_{\mathbf{A}}^{k}- R_h\overline{\psi}^{k-1}\overline{\theta}_{\mathbf{A}}^{k-1},\;\nabla\theta_{\psi}^{k-1}\right).}
\end{array}
\end{equation}

Using (\ref{eq:3-2-3}) and Young's inequality, we can estimate the first three terms on the right side of (\ref{eq:3-41}):
 \begin{equation}\label{eq:3-42}
 \begin{array}{@{}l@{}}
 {\displaystyle |\left(\nabla R_h\overline{\psi}^{m}\overline{\theta}_{\mathbf{A}}^{m},\;\theta_{\psi}^{m}\right)|
 +|\left(R_h\overline{\psi}^{m}\nabla\cdot\overline{\theta}_{\mathbf{A}}^{m},\;\theta_{\psi}^{m}\right)|
 +|\left(R_h\overline{\psi}^{0}\overline{\theta}_{\mathbf{A}}^{0},\;\nabla\theta_{\psi}^{0}\right)|}\\[2mm]
 {\displaystyle \quad\leq \|\nabla R_h\overline{\psi}^{m}\|_{\mathbf{L}^3}\|\overline{\theta}_{\mathbf{A}}^{m}\|_{\mathbf{L}^6}\|\theta_{\psi}^{m}\|_{\mathcal{L}^2}
 +\|I_{h}\overline{\psi}^{m}\|_{\mathcal{L}^{\infty}}\|\nabla\cdot\overline{\theta}_{\mathbf{A}}^{m}\|_{{L}^2}
 \|\theta_{\psi}^{m}\|_{\mathcal{L}^2}+C h^{2r} }\\[2mm]
{\displaystyle \quad\leq C\|\overline{\theta}_{\mathbf{A}}^{m}\|_{\mathbf{H}^1}\|\theta_{\psi}^{m}\|_{\mathcal{L}^2}
 +C h^{2r}\leq \frac{1}{32} D(\overline{\theta}_{\mathbf{A}}^{m},\overline{\theta}_{\mathbf{A}}^{m}) + C\|\theta_{\psi}^{m}\|^{2}_{\mathcal{L}^2}+C h^{2r}.}
 \end{array}
 \end{equation}

The last term on the right side of (\ref{eq:3-41}) can be estimated by
 \begin{equation}\label{eq:3-43}
 \begin{array}{@{}l@{}}
 {\displaystyle |\sum_{k=1}^{m}\left(R_h\overline{\psi}^{k}\overline{\theta}_{\mathbf{A}}^{k}- R_h\overline{\psi}^{k-1}\overline{\theta}_{\mathbf{A}}^{k-1},\;\nabla\theta_{\psi}^{k-1}\right)|}\\[2mm]
{\displaystyle \quad\leq  \Delta t \sum_{k=1}^{m}{\|\frac{1}{\Delta t}(R_h\overline{\psi}^{k}-R_h\overline{\psi}^{k-1})\|_{\mathcal{L}^3}\|\overline{\theta}_{\mathbf{A}}^{k}
\|_{\mathbf{L}^6}\|\nabla\theta_{\psi}^{k-1}\|_{\mathbf{L}^2}}}\\[2mm]
{\displaystyle \quad\quad+\Delta t\sum_{k=1}^{m}{\|R_h\overline{\psi}^{k-1}\|_{\mathcal{L}^{\infty}}
\|\partial \overline{\theta}_{\mathbf{A}}^{k}\|_{\mathbf{L}^2}\|\nabla\theta_{\psi}^{k-1}\|_{\mathbf{L}^2}}}\\[2mm]
{\displaystyle \quad\leq  C\Delta t \sum_{k=1}^{m}\big(\|\overline{\theta}_{\mathbf{A}}^{k}\|_{\mathbf{H}^{1}}
+\|\partial \overline{\theta}_{\mathbf{A}}^{k}\|_{\mathbf{L}^2}\big) \|\nabla\theta_{\psi}^{k-1}\|_{\mathbf{L}^2}}\\[2mm]
{\displaystyle \quad\leq  C\Delta t \sum_{k=0}^{m}\left(D({\theta}_{\mathbf{A}}^{k},{\theta}_{\mathbf{A}}^{k})
+\|\partial {\theta}_{\mathbf{A}}^{k}\|_{\mathbf{L}^2}^{2}+\|\nabla\theta_{\psi}^{k}\|_{\mathbf{L}^2}^{2}\right).}
 \end{array}
 \end{equation}

We thus get
 \begin{equation}\label{eq:3-44}
 \begin{array}{@{}l@{}}
 {\displaystyle |K_2|=|-\mathrm{i}\sum_{k=1}^{m}\left(R_h\overline{\psi}^{k}({\pi}_{h}\overline{\mathbf{A}}^{k}-\overline{\mathbf{A}}^{k}_{h}),\;
 \nabla\theta_{\psi}^{k}-\nabla\theta_{\psi}^{k-1}\right)|}\\[2mm]
 {\displaystyle \quad\leq \frac{1}{16} D(\overline{\theta}_{\mathbf{A}}^{m},\overline{\theta}_{\mathbf{A}}^{m}) + C\|\theta_{\psi}^{m}\|^{2}_{\mathcal{L}^2}+C h^{2r}}\\[2mm]
 {\displaystyle \quad+ C\Delta t \sum_{k=0}^{m}\left\{D({\theta}_{\mathbf{A}}^{k},{\theta}_{\mathbf{A}}^{k})
+\|\partial {\theta}_{\mathbf{A}}^{k}\|_{\mathbf{L}^2}^{2}+\|\nabla\theta_{\psi}^{k}\|_{\mathbf{L}^2}^{2}\right\}.}
 \end{array}
 \end{equation}

We similarly prove
 \begin{equation}\label{eq:3-45}
 \begin{array}{@{}l@{}}
 {\displaystyle |K_3|=|\mathrm{i}\sum_{k=1}^{m}\left(\nabla R_h\overline{\psi}^{k}({\pi}_{h}\overline{\mathbf{A}}^{k}
 -\overline{\mathbf{A}}^{k}_{h}),\;\theta_{\psi}^{k}-\theta_{\psi}^{k-1}\right)|}\\[2mm]
  {\displaystyle \quad\leq \frac{1}{16} D(\overline{\theta}_{\mathbf{A}}^{m},\overline{\theta}_{\mathbf{A}}^{m})+ C\|\theta_{\psi}^{m}\|^{2}_{\mathcal{L}^2}+C h^{2r}}\\[2mm]
   {\displaystyle \quad+ C\Delta t \sum_{k=0}^{m}\left(D({\theta}_{\mathbf{A}}^{k},{\theta}_{\mathbf{A}}^{k})
+\|\partial {\theta}_{\mathbf{A}}^{k}\|_{\mathbf{L}^2}^{2}+\|\nabla\theta_{\psi}^{k}\|_{\mathbf{L}^2}^{2}\right).}
 \end{array}
 \end{equation}

Adding (\ref{eq:3-40}), (\ref{eq:3-44}) and (\ref{eq:3-45}) together,  we find
\begin{equation}\label{eq:3-46}
\begin{array}{@{}l@{}}
{\displaystyle |\sum_{k=1}^{m} R_5^{k,3}| \leq  \frac{3}{16} D(\overline{\theta}_{\mathbf{A}}^{m},\overline{\theta}_{\mathbf{A}}^{m}) + C\|\theta_{\psi}^{m}\|^{2}_{\mathcal{L}^2}+C h^{2r}}\\[2mm]
 {\displaystyle \quad + C\Delta t \sum_{k=0}^{m}\left\{D({\theta}_{\mathbf{A}}^{k},{\theta}_{\mathbf{A}}^{k})
 +\|\partial {\theta}_{\mathbf{A}}^{k}\|_{\mathbf{L}^2}^{2}
 +\|\nabla\theta_{\psi}^{k}\|_{\mathbf{L}^2}^{2}\right\}.}
 \end{array}
 \end{equation}

It follows from (\ref{eq:3-30}), (\ref{eq:3-34}) and (\ref{eq:3-46}) that
 \begin{equation}\label{eq:3-47}
 \begin{array}{@{}l@{}}
 {\displaystyle \displaystyle |\Delta t\sum_{k=1}^{m}Q_{5}^{k}(\partial{\theta_{\psi}^{k}})|\leq C\big\{h^{2r}+(\Delta t)^{4}\big\}+\frac{3}{16} D(\overline{\theta}_{\mathbf{A}}^{m},\overline{\theta}_{\mathbf{A}}^{m}) + \frac{1}{16} \|\nabla\theta_{\psi}^{m}\|_{\mathbf{L}^2}^{2}
   }\\[2mm]
 {\displaystyle \qquad +C\|\theta_{\psi}^{m}\|_{\mathcal{L}^2}^{2} + C\Delta t \sum_{k=0}^{m}\left(D({\theta}_{\mathbf{A}}^{k},{\theta}_{\mathbf{A}}^{k}) +\|\partial {\theta}_{\mathbf{A}}^{k}\|_{\mathbf{L}^2}^{2}+\|\nabla\theta_{\psi}^{k}\|_{\mathbf{L}^2}^{2} \right)}\\[2mm]
 \end{array}
 \end{equation}

 Substituting (\ref{eq:3-19}), (\ref{eq:3-22}), (\ref{eq:3-23-0}), (\ref{eq:3-27}) and (\ref{eq:3-47}) into (\ref{eq:3-16-4}), we get
 \begin{equation}\label{eq:3-48}
\begin{array}{@{}l@{}}
{\displaystyle \frac{1}{2}B(\overline{\mathbf{A}}_{h}^{m};\theta_{\psi}^{m},\theta_{\psi}^{m}) \leq
 C \{h^{2r} +(\Delta t)^{4}\}+ \frac{3}{16} \|\nabla\theta_{\psi}^{m}\|_{\mathbf{L}^2}^{2}}\\[2mm]
 {\displaystyle \quad + \frac{1}{8}\Vert\nabla\overline{\theta}_{\phi}^{m}\Vert_{\mathbf{L}^2}^{2} +\frac{3}{16} D(\overline{\theta}_{\mathbf{A}}^{m},\overline{\theta}_{\mathbf{A}}^{m})
 +C\|\theta_{\psi}^{m}\|_{\mathcal{L}^2}^{2}}\\[2mm]
  {\displaystyle\quad +C \Delta t\sum_{k=0}^{m}\left\{ \|\nabla \theta_{\psi}^{k}\|^{2}_{\mathbf{L}^{2}} +\|\nabla \theta_{\phi}^{k}\|^{2}_{\mathbf{L}^{2}}+\|\partial \theta_{\phi}^{k}\|^{2}_{{L}^{2}}+ D(\theta_{\mathbf{A}}^{k},\theta_{\mathbf{A}}^{k}) + \|\partial \theta_{\mathbf{A}}^{k} \|^{2}_{\mathbf{L}^{2}}\right\}.}
 \end{array}
\end{equation}

Since
\begin{equation*}
B(\overline{\mathbf{A}}_{h}^{m};\theta_{\psi}^{m},\theta_{\psi}^{m}) = \|\nabla\theta_{\psi}^{m}\|_{\mathbf{L}^2}^{2} + \|\overline{\mathbf{A}}_{h}^{m}\theta_{\psi}^{m}\|_{\mathbf{L}^2}^{2} +
\left(f(\theta_{\psi}^{m},\theta_{\psi}^{m}),\;\overline{\mathbf{A}}_{h}^{m}\right),
\end{equation*}

 we have
 \begin{equation}\label{eq:3-56}
 \begin{array}{@{}l@{}}
 {\displaystyle \frac{5}{16}\|\nabla\theta_{\psi}^{m}\|_{\mathbf{L}^2}^{2}
 +\frac{1}{2}\|\overline{\mathbf{A}}_{h}^{m}\theta_{\psi}^{m}\|_{\mathbf{L}^2}^{2}
 \leq -\left(f(\theta_{\psi}^{m},\theta_{\psi}^{m}),\overline{\mathbf{A}}_{h}^{m}\right)
 +C\big\{h^{2r}+(\Delta t)^{4}\big\}}\\[2mm]
 {\displaystyle \quad + \frac{1}{8}\Vert\nabla\overline{\theta}_{\phi}^{m}\Vert_{\mathbf{L}^2}^{2} +\frac{3}{16} D(\overline{\theta}_{\mathbf{A}}^{m},\overline{\theta}_{\mathbf{A}}^{m})+C\|\theta_{\psi}^{m}\|_{\mathcal{L}^2}^{2} }\\[2mm]
 {\displaystyle \quad+C \Delta t\sum_{k=0}^{m}\left\{\|\nabla \theta_{\psi}^{k}\|^{2}_{\mathbf{L}^{2}} +\|\nabla \theta_{\phi}^{k}\|^{2}_{\mathbf{L}^{2}}+\|\partial \theta_{\phi}^{k}\|^{2}_{{L}^{2}}+ D(\theta_{\mathbf{A}}^{k},\theta_{\mathbf{A}}^{k}) + \|\partial \theta_{\mathbf{A}}^{k} \|^{2}_{\mathbf{L}^{2}}\right\}.}
 \end{array}
 \end{equation}

 It follows from (\ref{eq:3-84-0}) and the interpolation inequality (\ref{eq:3-2-1}) that
 \begin{equation}\label{eq:3-57}
 \begin{array}{@{}l@{}}
 {\displaystyle |\left(f(\theta_{\psi}^{m},\theta_{\psi}^{m}),\overline{\mathbf{A}}_{h}^{m}\right)|
 =|\frac{\mathrm{i}}{2}\left((\theta_{\psi}^{m})^{\ast}\nabla\theta_{\psi}^{m}-\theta_{\psi}^{m}\nabla (\theta_{\psi}^{m})^{\ast},\overline{\mathbf{A}}_{h}^{m}\right)|}\\[2mm]
 {\displaystyle \quad \leq \|\theta_{\psi}^{m}\|_{\mathcal{L}^3}
 \|\nabla\theta_{\psi}^{m}\|_{\mathbf{L}^2}
 \|\overline{\mathbf{A}}_{h}^{m}\|_{\mathbf{L}^6}
 \leq C\|\theta_{\psi}^{m}\|_{\mathcal{L}^3} \|\nabla\theta_{\psi}^{m}\|_{\mathbf{L}^2}}\\[2mm]
 {\displaystyle \quad \leq C\|\theta_{\psi}^{m}\|_{\mathcal{L}^3}^{2}
 +\frac{1}{32}\|\nabla\theta_{\psi}^{m}\|_{\mathbf{L}^2}^{2}
 \leq C\|\theta_{\psi}^{m}\|_{\mathcal{L}^2}
 \|\theta_{\psi}^{m}\|_{\mathcal{L}^6}+\frac{1}{32}\|\nabla\theta_{\psi}^{m}\|_{\mathbf{L}^2}^{2}}\\[2mm]
 {\displaystyle \quad \leq C\|\theta_{\psi}^{m}\|_{\mathcal{L}^2} \|\nabla\theta_{\psi}^{m}\|_{\mathbf{L}^2}
 +\frac{1}{32}\|\nabla\theta_{\psi}^{m}\|_{\mathbf{L}^2}^{2}}\\[2mm]
   {\displaystyle \quad
   \leq C\|\theta_{\psi}^{m}\|_{\mathcal{L}^2} ^{2}
   +\frac{1}{16}\|\nabla\theta_{\psi}^{m}\|_{\mathbf{L}^2}^{2}.}\\[2mm]
 \end{array}
 \end{equation}

We thus obtain
 \begin{equation}\label{eq:3-58}
  \begin{array}{@{}l@{}}
 {\displaystyle \frac{1}{4}\|\nabla\theta_{\psi}^{m}\|_{\mathbf{L}^2}^{2}
 \leq C\big\{h^{2r}+(\Delta t)^{4}\big\} + \frac{1}{8}\Vert\nabla\overline{\theta}_{\phi}^{m}\Vert_{\mathbf{L}^2}^{2} +\frac{3}{16} D(\overline{\theta}_{\mathbf{A}}^{m},\overline{\theta}_{\mathbf{A}}^{m})+C\|\theta_{\psi}^{m}\|_{\mathcal{L}^2}^{2}}\\[2mm]
 {\displaystyle \quad+C \Delta t\sum_{k=0}^{m}\left( \|\nabla \theta_{\psi}^{k}\|^{2}_{\mathbf{L}^{2}} +\|\nabla \theta_{\phi}^{k}\|^{2}_{\mathbf{L}^{2}}+\|\partial \theta_{\phi}^{k}\|^{2}_{{L}^{2}}+ D(\theta_{\mathbf{A}}^{k},\theta_{\mathbf{A}}^{k}) + \|\partial \theta_{\mathbf{A}}^{k} \|^{2}_{\mathbf{L}^{2}}\right).}
 \end{array}
 \end{equation}

 Multiplying (\ref{eq:3-15}) with $(C + 1)$ and adding to (\ref{eq:3-58}), we get
 \begin{equation}\label{eq:3-60}
 \begin{array}{@{}l@{}}
 {\displaystyle \frac{1}{4}\|\nabla\theta_{\psi}^{m}\|_{\mathbf{L}^2}^{2} +\|\theta_{\psi}^{m}\|_{\mathcal{L}^2}^{2}
 \leq C\big\{h^{2r}+(\Delta t)^{4}\big\} + \frac{1}{8}\Vert\nabla\overline{\theta}_{\phi}^{m}\Vert_{\mathbf{L}^2}^{2} +\frac{3}{16} D(\overline{\theta}_{\mathbf{A}}^{m},\overline{\theta}_{\mathbf{A}}^{m})}\\[2mm]
 {\displaystyle \quad+C \Delta t\sum_{k=0}^{m}\left( \|\nabla \theta_{\psi}^{k}\|^{2}_{\mathbf{L}^{2}} +\|\nabla \theta_{\phi}^{k}\|^{2}_{\mathbf{L}^{2}}+\|\partial \theta_{\phi}^{k}\|^{2}_{{L}^{2}}+ D(\theta_{\mathbf{A}}^{k},\theta_{\mathbf{A}}^{k}) + \|\partial \theta_{\mathbf{A}}^{k} \|^{2}_{\mathbf{L}^{2}}\right).}
 \end{array}
 \end{equation}

Adding (\ref{eq:3-84}) and (\ref{eq:3-92}) to (\ref{eq:3-60}), we end up with
\begin{equation}
\begin{array}{@{}l@{}}
{\displaystyle \|\theta_{\psi}^{m}\|_{\mathcal{L}^2}^{2}+\frac{1}{4}\|\nabla\theta_{\psi}^{m}\|_{\mathbf{L}^2}^{2} +\frac{1}{2}\|\partial \theta^{m}_{\phi}\|^{2}_{{L}^{2}} + \frac{1}{16}\|\nabla \theta^{m}_{\phi}\|^{2}_{\mathbf{L}^{2}}}\\[2mm]
{\displaystyle\quad +\frac{1}{2}\|\partial \theta_{\mathbf{A}}^{m}\|_{\mathbf{L}^2}^{2} + \frac{1}{8}D({\theta}_{\mathbf{A}}^{m}, {\theta}_{\mathbf{A}}^{m})
\leq C\big\{h^{2r}+(\Delta t)^{4}\big\}}\\[2mm]
{\displaystyle \quad  + C \Delta t\sum_{k=1}^{m}\big\{ \|\nabla \theta_{\psi}^{k}\|^{2}_{\mathbf{L}^{2}} +\|\nabla \theta_{\phi}^{k}\|^{2}_{\mathbf{L}^{2}}+\|\partial \theta_{\phi}^{k}\|^{2}_{{L}^{2}}}\\[2mm]
{\displaystyle \quad + D(\theta_{\mathbf{A}}^{k},\theta_{\mathbf{A}}^{k}) + \|\partial \theta_{\mathbf{A}}^{k} \|^{2}_{\mathbf{L}^{2}}\big\}.}
\end{array}
\end{equation}

Now by applying the discrete Gronwall's inequality and choosing a sufficiently small $ \Delta t $ such that $C\Delta t \leq \frac{\displaystyle 1}{\displaystyle 2}$, we conclude
\begin{equation}
\begin{array}{@{}l@{}}
{\displaystyle \max_{0 \leq k \leq m}\left\{\|\theta_{\psi}^{m}\|_{\mathcal{L}^2}^{2}+\|\nabla\theta_{\psi}^{m}\|_{\mathbf{L}^2}^{2} +\|\partial \theta^{m}_{\phi}\|^{2}_{{L}^{2}} + \|\nabla \theta^{m}_{\phi}\|^{2}_{\mathbf{L}^{2}} +\|\partial \theta_{\mathbf{A}}^{m}\|_{\mathbf{L}^2}^{2} + D({\theta}_{\mathbf{A}}^{m}, {\theta}_{\mathbf{A}}^{m})\right\}} \\[2mm]
{\displaystyle \qquad \qquad \leq  C\exp\left(\frac{TC}{1-C \Delta t}\right) \left\{h^{2r}+(\Delta t)^{4}\right\} \leq C\exp\left(2TC\right) \left\{h^{2r}+(\Delta t)^{4}\right\} .}
\end{array}
\end{equation}

If we take $C_{\ast} \geq  C\exp\left(2TC\right)$, then (\ref{eq:3-2-6}) holds for $k = m$. By virtue of the mathematical induction,
we complete the proof of (\ref{eq:3-2-6}). Furthermore, using the triangle inequality and the equality $ \theta_{\mathbf{A}}^{M}=\Delta t \sum\limits_{k=1}^{M}{\partial \theta_{\mathbf{A}}^{k}}+ \theta_{\mathbf{A}}^{0} $,
we can complete the proof of Theorem~\ref{thm2-1}.

\section{Numerical testes}\label{sec-4}
In this section, we present numerical experiments to illustrate the error estimates.

We consider the following Maxwell--Schr\"{o}dinger's equations:
\begin{equation}\label{eq:4-1}
\left\{
\begin{array}{@{}l@{}}
{\displaystyle  -\mathrm{i}\frac{\partial \psi}{\partial t}+
\frac{1}{2}\left(\mathrm{i}\nabla +\mathbf{A}\right)^{2}\psi
 + V_{0}\psi + \phi\psi= f(\mathbf{x},t) ,\,\, (\mathbf{x},t)\in
\Omega\times(0,T),}\\[2mm]
{\displaystyle \frac{\partial ^{2}\mathbf{A}}{\partial t^{2}}+\nabla\times
(\nabla\times \mathbf{A}) - \nabla(\nabla \cdot \mathbf{A})+\frac{\mathrm{i}}{2}\big(\psi^{*}\nabla{\psi}-\psi\nabla{\psi}^{*}\big) }\\[2mm]
{\displaystyle \qquad\qquad+\vert\psi\vert^{2}\mathbf{A}=\mathbf{g}(\mathbf{x},t),
\,\,\quad (\mathbf{x},t)\in \Omega\times(0,T),}\\[2mm]
{\displaystyle  \frac{\partial ^{2}\phi}{\partial t^{2}}-\Delta \phi - \vert\psi\vert^{2} = l(\mathbf{x},t),\,\, (\mathbf{x},t)\in \Omega\times(0,T) }.
\end{array}
\right.
\end{equation}
where the initial-boundary conditions are given in (\ref{eq:1-3})-(\ref{eq:1-4}).

Let $\Omega = (0,1)^3$, $T = 4$ and $V_0 = 5$. The exact solution $(\psi,\mathbf{A},\phi)$ of (\ref{eq:4-1}) is defined by
\begin{equation*}
\begin{array}{@{}l@{}}
{\displaystyle \psi(\mathbf{x}, t) =  (1+0.5 t)e^{\mathrm{i}\pi t}  \sin(2\pi x_1)\sin(2\pi x_2)\sin(2\pi x_3), }
\end{array}
\end{equation*}
\begin{equation*}
\begin{array}{@{}l@{}}
{\displaystyle \mathbf{A}(\mathbf{x},t)= \cos(\pi t)\Big(\cos(\pi x_1)\sin(\pi x_2)\sin(\pi x_3),}\\[2mm]
{\displaystyle \sin(\pi x_1)\cos(\pi x_2)\sin(\pi x_3), \;\;\sin(\pi x_1)\sin(\pi x_2)\cos(\pi x_3)\Big),}
\end{array}
\end{equation*}
\begin{equation*}
\phi(\mathbf{x}, t) = \big(t+sin(\pi t)\big)x_1x_2x_3(1.0-x_1)(1.0-x_2)(1.0-x_3).
\end{equation*}

The functions  $f(\mathbf{x},t)$, $\mathbf{g}(\mathbf{x},t)$ and $l(\mathbf{x},t)$ at the right hand side of  (\ref{eq:4-1}) are chosen correspondingly to the exact solution $(\psi,\mathbf{A},\phi)$.

We partition a whole domain $ \Omega $ into quasi-uniform tetrahedrons with $M +1$ nodes in each direction and $6M^{3}$ elements in total.
The system (\ref{eq:4-1}) is solved by the proposed Crank-Nicolson Galerkin finite element scheme (\ref{eq:2-9}) with linear elements and quadratic elements, respectively. 
To confirm the convergence rate of the proposed method, we take $\Delta t = h^{\frac{1}{2}}$ for the linear element method and $\Delta t = h$ for the quadratic element method respectively. 
The numerical results for the linear element method and the quadratic element method at time $t=1.0, 2.0, 3.0, 4.0 $  are displayed in Tables~\ref{table4-1} and~\ref{table4-2}, respectively. 
It clearly shows that they are in good agreement with the error estimates presented in Theorem~\ref{thm2-1}. 

\begin{table}[htb]
\caption{$H^1$ error of linear FEM with $h = \frac{1}{M}$ and $\Delta t = h^{\frac12}$.}\label{table4-1}
\begin{center}
\begin{tabular}{c|cccc}
   \hline\hline
    & &   $\Vert \mathbf{A}^{k}_{h} - \mathbf{A}(\cdot, t_k) \Vert _{\mathbf{H}^{1}} $ \\
   \hline
    t &  M=25 & M=50 & M=100 & Order   \\
    \hline
    1.0 & 4.9855e-01 & 2.3894e-01 & 1.1887e-01 & 1.03 \\
    \hline
     2.0 & 6.7234e-01 & 3.2057e-01 & 1.7574e-01 & 0.97 \\
     \hline
     3.0 & 4.6375e-01 & 2.2119e-01 & 1.0373e-01 & 1.08 \\
    \hline
     4.0 & 5.8205e-01 & 3.0487e-01 & 1.4287e-01 & 1.02 \\
     \hline\hline
   & &   $ \Vert \psi_{h}^{k}-\psi(\cdot,t_k)\Vert_{\mathcal{H}^1}$    \\
     \hline
    t &  M=25 & M=50 & M=100 & Order   \\
    \hline
    1.0 & 3.0718e-01 & 1.4419e-01 & 7.4104e-02 & 1.02 \\
    \hline
     2.0 & 4.3713e-01 & 2.1289e-01 & 1.1430e-01 & 0.97 \\
     \hline
     3.0 & 3.0004e-01 & 1.4202e-01 & 6.9620e-02 & 1.05 \\
    \hline
     4.0 & 2.2543e-01 & 1.2316e-01 & 6.0316e-02  & 0.95 \\
     \hline\hline
   & &   $ \Vert \phi_{h}^{k}-\phi(\cdot,t_k)\Vert_{{H}^1}$    \\
     \hline
    t &  M=25 & M=50 & M=100 & Order   \\
    \hline
    1.0 & 1.9412e-01 & 8.6435e-02 & 4.3191e-02 & 1.07 \\
    \hline
     2.0 & 1.2473e-01 & 7.0233e-02 & 3.3148e-02 & 0.96 \\
     \hline
     3.0 & 1.1031e-01 & 6.1394e-02 & 2.9217e-02 & 0.96 \\
    \hline
     4.0 & 8.7695e-02 & 4.1380e-02 & 2.1815e-02  & 1.00 \\
     \hline\hline
  \end{tabular}
  \end{center}
\end{table}

\begin{table}[htb]
\caption{$H^1$ error of quadratic FEM with $h = \Delta t= \frac{1}{M}.$}\label{table4-2}
\begin{center}
\begin{tabular}{c|cccc}
   \hline\hline
    & &   $\Vert \mathbf{A}^{k}_{h} - \mathbf{A}(\cdot, t_k) \Vert _{\mathbf{H}^{1}} $ \\
   \hline
    t &  M=25 & M=50 & M=100 & Order   \\
    \hline
    1.0 & 6.5062e-02   & 1.6679e-02 & 4.1805e-03 & 1.98 \\
    \hline
     2.0 & 5.5496e-02 & 1.5018e-02  & 3.7016e-03 & 1.95 \\
     \hline
     3.0 & 4.1903e-02 & 1.3291e-02 & 2.4825e-03 & 2.04 \\
    \hline
     4.0 & 3.4765e-02 &  8.4562e-03 & 2.0952e-03  & 2.03 \\
     \hline\hline
   & &   $ \Vert \psi_{h}^{k}-\psi(\cdot,t_k)\Vert_{\mathcal{H}^1}$    \\
     \hline
    t &  M=25 & M=50 & M=100 & Order   \\
    \hline
    1.0 & 1.1930e-02 & 3.2391e-03 & 8.0967e-04 & 1.94 \\
    \hline
     2.0 & 1.0237e-02 & 2.9786e-03 & 5.9668e-04 & 2.05 \\
     \hline
     3.0 & 2.2342e-02 & 4.8413e-03 & 1.3683e-03 & 2.01 \\
    \hline
     4.0 & 1.3203e-02 & 3.1681e-03 & 7.8548e-04  & 2.04 \\
     \hline\hline
     & &   $ \Vert \phi_{h}^{k}-\phi(\cdot,t_k)\Vert_{{H}^1}$    \\
     \hline
    t &  M=25 & M=50 & M=100 & Order   \\
    \hline
    1.0 & 3.1070e-02 & 7.8114e-03 & 1.8648e-03 & 2.03 \\
    \hline
     2.0 & 2.7189e-02 & 8.2275e-03 & 1.8301e-03 & 1.95 \\
     \hline
     3.0 & 3.1162e-02 & 6.9836e-03 & 1.7682e-03 & 2.07 \\
    \hline
     4.0 & 2.5343e-02 & 6.4257e-03 & 1.6146e-03  & 1.99 \\
     \hline\hline
  \end{tabular}
  \end{center}
\end{table}

\clearpage

\end{document}